\documentclass[12pt]{amsart}
\usepackage[all]{xy} 
\usepackage{amsmath}
\usepackage{amsfonts}
\usepackage{amssymb}

\DeclareFontEncoding{OT2}{}{} 
\newcommand{\textcyr}[1]{%
 {\fontencoding{OT2}\fontfamily{wncyr}\fontseries{m}\fontshape{n}
 \selectfont #1}}
\newcommand{\Sha}{{\mbox{\textcyr{Sh}}}}

\def\ut{\widetilde{U}}

\def\bz{{\mathbb Z}\,}
\def\bh{{\mathbb H}}
\def\bq{{\mathbb Q}}
\def\bg{{\mathbb G}}
\def\spec{{\rm{Spec}}\,}
\def\pdiv{\text{$p\kern 0.1em$-div}}

\def\be{\kern -.1em}
\def\lbe{\kern -.05em}
\def\s{\mathcal }
\def\ra{\rightarrow}
\def\e{\kern 0.08em}
\def\le{\kern 0.04em}
\def\ng{\kern -0.04em}

\def\krn{{\rm{Ker}}\e }
\def\img{{\rm{Im}}\e }

\newtheorem{lemma}{Lemma}[section]
\newtheorem{theorem}[lemma]{Theorem}
\newtheorem{corollary}[lemma]{Corollary}
\newtheorem{proposition}[lemma]{Proposition}
\theoremstyle{definition}

\theoremstyle{remark}
\newtheorem{remark}[lemma]{Remark}
\newtheorem{remarks}[lemma]{Remarks}

\begin{document}

\title[Duality for 1-motives over function fields]{Arithmetic duality
theorems for 1-motives over function fields}


\author{Cristian D. Gonz\'alez-Avil\'es}
\address{Departamento de Matem\'aticas, Universidad de La Serena,
La Serena, Chile} \email{cgonzalez@userena.cl}


\thanks{The author was partially supported by Fondecyt grant
1061209}

\maketitle

\begin{abstract} In this paper we obtain a Poitou-Tate exact sequence
for finite and flat group schemes over a global function field. In
addition, we extend the duality theorems for 1-motives over number
fields obtained by D.Harari and T.Szamuely to the function field
case.
\end{abstract}

\section{Introduction}
Let $X$ be a smooth projective curve over a finite field of
characteristic $p$ and let $K$ be the function field of $X$. In this
paper we establish a Poitou-Tate exact sequence for finite and flat
group schemes of $p$-power order over $K$, thereby extending the
well-known Poitou-Tate exact sequence in Galois cohomology [15,
Theorem I.4.10, p.70]. See Theorem 4.12 below for the precise
statement. In particular, we obtain the following duality theorem.

\begin{theorem} Let $N$ be a $p$-primary finite and flat group scheme
over $K$. Then there exists a perfect pairing of finite groups
$$
\Sha^{\e 1}(K,N)\times\!\! \Sha^{\e 2}\!\be\left(K,N^{\e d}\e\right)
\ra\bq_{p}/\e\bz_{\be p}\,.
$$
\end{theorem}
The Tate-Shafarevich groups appearing in the statement of the
theorem are defined in terms of flat cohomology.

In addition, we extend the duality theorems for 1-motives over
number fields obtained by D.Harari and T.Szamuely [10] to the
function field case. In particular, we obtain the following result.

\begin{theorem} Let $M$ be a 1-motive over $K$ with dual 1-motive
$M^{*}$. Then there exists a canonical pairing
$$
\Sha^{\e 1}(K,M)(p)\e\times\!\!\Sha^{\e 1}(K,M^{*})(p)\ra\bq/\e\bz
$$
whose left and right kernels are the maximal divisible subgroups of
each group.
\end{theorem}

The paper has 6 Sections. In Section 2 we prove some elementary
results needed in the sequel. Section 3 is a brief summary of the
facts that we need on 1-motives (readers wishing to learn more about
the theory of 1-motives are advised to read [1]). In Section 4,
which is independent of Sections 5 and 6, we establish the
Poitou-Tate exact sequence for $p$-primary finite and flat group
schemes. In Section 5 we prove an ``integral version" of Theorem
1.2, namely an analogous statement with $\spec\e K$ replaced by an
open affine subset $U$ of $X$. Theorem 1.2 is then deduced from this
integral version in Section 6 by passing to the limit as $U$ shrinks
to $\spec K$.

The methods of this paper yield both a general Poitou-Tate exact
sequence and a Cassels-Tate dual exact sequence for 1-motives over
global fields (extending the results of [10, \S 5], [11, \S 5] and
[7]). These sequences require a significant amount of extra work in
relation to [op.cit.] and will be established in separate
publications.

\section*{Acknowledgements} I am very grateful to D. Harari for
making several helpful comments, carefully reading a preliminary
version of this paper and pointing out an error in it. I also thank
S. J. Edixhoven, R. de Jong, J. S. Milne, L. Taelman and the
anonymous referee for additional helpful comments, corrections and
valuable suggestions. Most of this paper was written during my
visits to Universit\'e de Bordeaux I, ICTP Trieste and Universiteit
Leiden in the summer of 2007 as an Erasmus Mundus ALGANT Scholar. I
am very grateful to S. J. Edixhoven for the invitation to
participate in the ALGANT program, to the institutions named above
for their warm hospitality and to the European Union for financial
support.

\section{Preliminaries}

Let $K$ be a global function field of characteristic $p>0$. For
basic information on these fields, and their completions, the reader
is referred to [4, Chapters I and II]. For any prime $v$ of $K$,
$K_{v}$, ${\s O}_{v}$ and $k(v)$ will denote, respectively, the
completion of $K$ at $v$, its ring of integers and the corresponding
residue field. Thus ${\s O}_{v}$ is a complete discrete valuation
ring. We will write $X$ for the unique smooth complete curve over
the field of constants of $K$ having $K$ as its function field. The
primes of $K$ will often be identified with the closed points of
$X$. A direct product extending over all $v\in U$ for some nonempty
open subset $U$ of $X$ is to be understood as extending over all
{\it closed} points of $U$.

For any abelian group $B$ and positive integer $n$, we will write
$B_{n}$ for the $n$-torsion subgroup of $B$ and $B/n$ for the
quotient $B/nB$. Further, we will write $B(p)=\cup_{\e m\geq
1}B_{p^{\le m}}$ (the $p$-primary torsion subgroup of $B$),
$B^{(p)}=\varprojlim_{m}B/p^{\le m}$ (the $p$-adic completion of
$B$) and $T_{p}\e B=\varprojlim_{m}B_{\e p^{\le m}}$ (the $p$-adic
Tate module of $B$). Also, we set $B_{\e p\e\text{-div}}=\bigcap_{\,
m}\e p^{m}B$ (the subgroup of $B$ of infinitely $p$-divisible
elements). If $B$ is a $p$-primary torsion abelian group of finite
cotype (i.e., if $B$ is isomorphic to a direct sum
$(\bq_{p}/\e\bz_{\be p})^{\le r}\oplus F$ for some integer $r\geq 0$
and some finite group $F$ or, equivalently, if $B_{\e p}$ is
finite), then $B_{\e p\e\text{-div}}$ coincides with the maximal
$p$-divisible subgroup of $B$ (and $B/B_{\e p\e\text{-div}}$ is
finite). For simplicity, we will often write $B/\text{$p\e$-div}$
for $B/B_{p\e\text{-div}}$.

\begin{lemma} Let $B$ be any abelian group. Then the canonical map
$B\ra B^{\e(p)}$ induces an injection
$$
B(p)/\text{$p\kern 0.1em$-{\rm{div}}}\hookrightarrow B^{\e(p)}(p).
$$
\end{lemma}
\begin{proof} There exists an exact sequence of inverse systems
$0\ra p^{\e m}B\ra  B\ra B/p^{\e m}\ra 0$, where the transition maps
are $p^{\e m+1}B\hookrightarrow p^{\e m}B$, $\text{Id}\colon B\ra B$
and $\text{proj}\colon B/p^{\e m+1}\ra B/p^{\e m}$. Taking inverse
limits, we obtain an exact sequence
$$
0\ra B_{\e p\e\text{-div}}\ra B\ra B^{\e(p)}
$$
and therefore an exact sequence
$$
0\ra B_{\e p\e\text{-div}}(p)\ra B(p)\ra B^{\e(p)}(p).
$$
Since $B_{\e p\e\text{-div}}(p)=B(p)_{p\e\text{-div}}$, the lemma
follows.
\end{proof}

If $B$ is any abelian topological group, we will write $\widehat{B}$
(or $B\,\widehat{\,\phantom{.}}\,$) for $\varprojlim_{\e I\in\s
S}B/I$, where $\s S$ is the family of open subgroups of $B$ of
finite {\it $p$-power} index. Clearly, if $B$ is a discrete torsion
abelian group of finite cotype, then $\widehat{B}=B/\pdiv$. There
exists a canonical isomorphism
$\widehat{B}=\big(B^{(p)}\big)\widehat{\,\phantom{.}}\,$, whence
there exists a canonical map $B^{(p)}\ra \widehat{B}$. We set
$B^{D}=\text{Hom}_{\e\text{cont.}}(B,\bq_{p}/\e\bz_{\be p})$, where
$\bq_{p}/\e\bz_{\be p}$ is endowed with the discrete topology. We
endow $B^{D}$ with the compact-open topology, i.e., the open subsets
of $B^{D}$ are arbitrary unions of finite intersections of sets of
the form $\{f\in B^{D}\colon f(K)\subset U\}$, where $K\subset B$ is
compact and $U\subset\bq_{p}/\e\bz_{\be p}$ is open (i.e.,
arbitrary). Note that, if $B$ is discrete and finitely generated,
then $B^{D}=\left(B^{(p)}\right)^{\be D}=\widehat{B}^{\e D}$ is a
discrete $p$-primary torsion group and $B^{DD}=B^{(p)}=\widehat{B}$.

\smallskip

A pairing of discrete abelian groups $A\times B\ra\bq_{p}/\e\bz_{\be
p}$ is called {\it non-degenerate on the right} (resp. {\it left})
if the induced homomorphism $B\ra A^{D}$ (resp. $A\ra B^{D}$) is
injective. It is called {\it non-degenerate} if it is non-degenerate
both on the right and on the left. The pairing is said to be {\it
perfect} if the homomorphisms $B\ra A^{D}$ and $A\ra B^{D}$ are
isomorphisms. It is not difficult to see that a perfect pairing
$A\times B\ra\bq_{p}/\e\bz_{\be p}$ induces pairings
$A(p)\times(B/p\e\text{-div})\ra\bq_{p}/\e\bz_{\be p}$ and
$(A/p\e\text{-div})\times B(p)\ra \bq_{p}/\e\bz_{\be p}$ which are
non-degenerate on the left and on the right, respectively.

\smallskip

\begin{lemma} Let $p$ be a prime number.

\begin{enumerate}

\item[(a)] Let $B$ be an abelian group and let $A$ be a torsion
subgroup of
$B$. If $B_{p}=0$, then $(B/A)_{p}=0$.
\item[(b)] Let $A\ra B\ra C\ra 0$ be an exact sequence of discrete
torsion abelian groups. Then the induced sequence $0\ra C(p)^{D}\ra
B(p)^{D}\ra A(p)^{D}$ is exact.
\end{enumerate}
\end{lemma}
\begin{proof} Since every element of $A$ is annihilated by an
integer which is prime to $p\,$, $A$ is $p$-divisible. This implies
(a). Assertion (b) follows from the fact that $B\ra B^{D}$ is an
exact functor on the category of discrete abelian groups.
\end{proof}

\smallskip

In this paper we consider only {\it commutative} group schemes,
and therefore the qualification ``commutative" will often be
omitted when discussing group schemes. Further, all cohomology
groups below are flat (fppf) cohomology groups.

\smallskip

Now let $N$ be a finite, flat (commutative) group scheme over $\spec
K$ and let $\s F=\s F(N)$ be the set of all pairs $(U,\s N\e)$,
where $U$ is a nonempty open {\it affine} subscheme of $X$ (i.e.,
$U\neq X$) and $\s N$ is a finite and flat group scheme over $U$
which extends $N$, i.e., $\s N\times_{U}\e\spec K=N$. Then $\s F$ is
a nonempty [13, p.294] directed and partially ordered
set\footnote{$\e$Note that, if $(U,\s N\e),(V,{\s
N}^{\e\prime}\e)\in\s F$, then there exists a nonempty open subset
of $U\cap V$ over which $\s N$ and ${\s N}^{\e\prime}$ are
isomorphic.} with the partial ordering $(U,\s N\e)\leq (V,{\s
N}^{\e\prime}\e)$ if and only if $V\subset U$ and ${\s N}\!\mid_{\e
V}={\s N}^{\e\prime}$. Clearly, $\varprojlim_{(U,\e\s N\e)\in\s F}\e
U=\bigcap_{\e(U,\e\s N\e)\in\s F}\e U=\spec K$.

\begin{lemma} With the above notations, for every $i\geq 0$ the
canonical map
$$
\varinjlim_{\e(U,\e\s N\e)\in\s F}H^{\e i}(U,\s N\e)\ra H^{\e i}(K,N)
$$
is an isomorphism.
\end{lemma}
\begin{proof} The result is clear if $i=0$. Assume now that $i\geq 1$. If $(U,\s N\e)\in\s F$, then $\s N$ admits a {\it canonical} resolution
$$
0\ra\s N\ra \s G_{\e 0}\ra \s G_{\e 1}\ra 0,
$$
where $\s G_{\e 0}$ and $\s G_{\e 1}$ are smooth affine group
schemes of finite type over $U$. See [2, \S 2.2.1, p.25]. Now, using
the fact that flat and \'etale cohomology coincide on smooth group
schemes [14, Theorem 3.9, p.114], we obtain the following exact
sequence which is functorial in $(U,\s N\e)$:
$$
H^{ i-1}_{\text{\'et}}(U,\s G_{\e 0}\e)\ra H^{ i-1}_{\text{\'et}}
(U,\s G_{\e 1}\e)\ra H^{\e i}(U,\s N\e)\ra H^{ i}_{\text{\'et}}(U,\s
G_{\e 0}\e)\ra H^{ i}_{\text{\'et}} (U,\s G_{\e 1}\e)
$$
An analogous exact sequence exists over $K$, and these exact
sequences form the top and bottom row, respectively, of a natural
exact commutative diagram. Since the canonical maps
$$
\varinjlim_{(U,\e\s N\e)\in\s F}\e H^{\e j}_{\text{\'et}}(U,
\s G_{\e l})\ra
H^{\e j}_{\text{\'et}}(K,G_{l})
$$
are isomorphisms for $j=i-1$ or $i$ and $l=0$ or 1 by [8, Theorem
VII.5.7, p.361], the five-lemma applied to the direct limit of the
diagram mentioned above yields the desired result.
\end{proof}

\begin{lemma} There exists a canonical isomorphism
$$
H^{\e 1}(K,N)^{D}=\displaystyle\varprojlim_{(U,\e\s N\e)\in\s
F}H^{\e 1}(U,\s N\e)^{D}.
$$
\end{lemma}
\begin{proof} Let $(U,\e\s N\e)\in\s
F$ and let $V$ be a nonempty open subset of $U$. Since $H^{\e
1}_{v}(O_{v},\s N\e )=0$ for any $v$ [15, beginning of \S III.7,
p.349], the localization sequence for the pair $V\subset U$
[op.cit., Proposition III.0.3(c), p.270] shows that the canonical
map of discrete groups $H^{\e 1}(U,\s N\e)\ra H^{\e 1}(V,\s N\e)$ is
injective. The lemma now follows from Lemma 2.3 above and [16,
Propositions 3 and 7].
\end{proof}

\smallskip

Let $S$ be a scheme. An {\it $S$-torus} $\s T$ is a smooth $S$-group
scheme which, locally for the \'etale topology on $S$, is isomorphic
to ${\bg}_{m}^{r}$ for some positive integer $r$.

We will use without explicit mention the fact that the \'etale and
fppf cohomology groups of a smooth group scheme coincide. In
particular, $H^{i}(S,\s T)=H^{i}_{\text{\'et}}(S,\s T)$ for an
$S$-torus $\s T$ as above.

Finally, for each prime $v$ of $K$, $\text{inv}_{\lbe v}\,\colon
\text{Br}(K_{v})\ra\bq/\e\bz$ will denote the usual invariant map of
local class field theory.

\section{Generalities on 1-motives}

Let $S$ be a scheme. We will write $S_{\e\text{fppf}}$ for the small
fppf site over $S$, ${\s F}_{S}$ for the category of abelian sheaves
on $S_{\e\text{fppf}}$, ${\s C}^{b}({\s F}_{S})$ for the category of
bounded complexes of objects in ${\s F}_{S}$ and ${\s D}^{\e b}( {\s
F}_{S})$ for the associated derived category.

Recall that a (smooth) 1-{\it motive}
$M=(Y,A,T,G,u)$ over $S$ consists of the following data:

\medskip

\begin{enumerate}

\item[1.] An $S$-group scheme $Y$ which, locally for the \'etale
topology on $S$, is isomorphic to ${\bz}^{\be r}$ for some $r\geq
0$.

\medskip

\item[2.] A commutative $S$-group scheme $G$ which is an extension
of an abelian $S$-scheme $A$ by an $S$-torus $T$:
$$
0\ra T\ra G\overset{\pi}\longrightarrow A\ra 0.
$$

\medskip

\item[3.] An $S$-homomorphism $u\colon Y\ra G$.

\end{enumerate}

\bigskip

We will often identify $M$ with the mapping cone of $u$, i.e.,
$M=C^{\bullet}(u)=(Y\overset{u}\longrightarrow G\e)$, with $Y$ placed
in degree $-1$ and $G$ placed in degree 0. Thus there exists a
distinguished triangle
\begin{equation}
Y\ra G\ra M\ra Y[1].
\end{equation}
Note also that $M$ defines, in a canonical way, an object of $D^{\e
b}(\text{fppf})$. Every 1-motive $M$ comes equipped with a natural
increasing 3-term {\it weight filtration}: $W_{i}(M)=0$ for $i\leq
-3$, $W_{-2}(M)=(0\ra T)$, $W_{-1}(M)=(0\ra G\e)$ and $W_{i}(M)=M$
for $i\geq 0$. The 1-motive
$$
M^{\prime}=M/W_{-2}(M)=(Y\overset{h}\longrightarrow A),
$$
where $h=\pi\circ u$, will play an auxiliary role below. It fits
into an exact sequence
\begin{equation}
0\ra T\ra M\ra M^{\e\prime}\ra 0,
\end{equation}
where $T$ is regarded as a complex concentrated in degree zero. Now,
to each 1-motive $M=(Y,A,T,G,u)$ as above, one can associate its
Cartier dual $M^{*}=(Y^{*},A^{*},T^{*},G^{*},u^{*})$. Here $Y^{*}$
is the sheaf of characters of $T$, $A^{*}$ is the abelian scheme
dual to $A$ and $T^{*}$ is the $S$-torus with group of characters
$Y$. The $S$-group scheme $G^{*}$ associated to $M$ may be
constructed as follows. Assume first that $M=M^{\prime}$ (i.e.,
$T=0$). In this case $M^{*}=(M^{\e\prime})^{*}=(0\ra G^{*})$, where
$G^{*}$ is the $S$-group scheme which represents the functor
$S^{\e\prime}\mapsto\text{Ext}^{1}_{S^{\prime}}(M^{\e\prime},{\mathbb
G}_{m})$ on ${\s C}^{b}({\s F}_{S})$ (the representability of this
functor follows from the generalized Weil-Barsotti
formula\footnote{The referee calls attention to the fact that the
Weil-Barsotti formula over imperfect fields has not been properly
established in the literature. However, O.Wittenberg has recently
given the necessary missing details in the Appendix to [18].}. The
1-motive $(M^{\prime})^{*}$ is naturally endowed with a biextension
(in the sense of [5, 10.2.1, p.60]) ${\s P}^{\e\prime}$ of
$(M^{\e\prime},(M^{\e\prime})^{*})$ by ${\mathbb G}_{m}$, namely the
pullback of the canonical Poincar\'e biextension of $(A,A^{*})$ by
${\mathbb G}_{m}$ under the map $f^{\e\prime}\times g^{\e\prime}$,
where $f^{\e\prime}=(0,\text{Id})\colon M^{\e\prime}=(Y\ra A)\ra A$
and $g^{\e\prime}$ is the composite
$(M^{\e\prime})^{*}\overset{(0,\text{Id})}\longrightarrow
G^{*}\overset{\pi^{*}}\longrightarrow A^{*}$. Now let $M$ be an
arbitrary 1-motive. By (2), $M$ represents a class in
$\text{Ext}^{1}_{S}(M^{\e\prime},T)$. Thus any
$\chi_{_{S^{\e\prime}}}\in
Y^{*}(S^{\e\prime})=\text{Hom}_{S^{\e\prime}}(T,{\mathbb G}_{m})$
induces an element
$u^{*}(\chi_{_{S^{\e\prime}}})=(\chi_{_{S^{\e\prime}}})_{*}
(M_{\be_{S^{\e\prime}}})\in\text{Ext}^{1}_{S^{\e\prime}}(
M^{\e\prime},{\mathbb G}_{m})=G^{*}(S^{\e\prime})$, which defines an
$S$-homomorphism $u^{*}\colon Y^{*}\ra G^{*}$. The associated
1-motive $M^{*}=(Y^{*}\overset{u^{*}}\longrightarrow G^{*})$ is the
{\it Cartier dual} of $M$. The corresponding biextension ${\s P}$ of
$(M,M^{*})$ by ${\mathbb G}_{m}$ is the pullback of ${\s P}^{
\e\prime}$ under the map $f\times g$, where $f=(\text{Id},\pi)\colon
M=(Y\ra G\e)\ra M^{\e\prime}=(Y\ra A\e)$ and $g=(0,\text{Id})\colon
M^{*}=(Y^{*}\ra G^{*})\ra(M^{\e\prime})^{*}=(0\ra G^{*})$ are the
natural maps.

Now, as in [9, VII.3.6.5], (the isomorphism class of) $\s P$
corresponds to a map $M\otimes^{\,\mathbf{L}}\! M^{*}\ra {\mathbb
G}_{m}[1]$ in ${\s D}^{\e b}( {\s F}_{S})$. This map in turn induces
pairings
$$
{\bh}^{\e i}(S,M)\times {\bh}^{\e j}(S,M^{*})\ra {\bh}^{\e
i+j+1}(S,{\mathbb G}_{m})
$$
for each $i,j\geq -1$.

Next, for any positive integer $n$, let
\begin{equation}
T_{\bz/n}(M)={\bh}^{-1}(C^{\e\bullet}(n))={\bh}^{\e
0}(M[-1]\!\otimes^{\mathbf{L}} \bz/n),
\end{equation}
where $C^{\e\bullet}(n)$ is the mapping cone of the
multiplication-by-$n$ map on $M$ (to verify the second equality in
(3), use the fact that $\bz/n$ is quasi-isomorphic to the complex of
flat modules $(\bz\overset{n}\longrightarrow\bz)$). It is a finite
and flat $S$-group scheme which fits into an exact sequence
$$
0\ra G_{n}\ra T_{\bz/n}(M)\ra Y/n\ra 0.
$$
See [1, \S 2.3, p.12]. It is not difficult to see that
$T_{\bz/n}(M)$ is the sheaf associated to the presheaf
$S^{\e\prime}\mapsto {\s F}_{\,\bz/n}(M)(S^{\e\prime})$, where
\begin{equation}
{\s F}_{\,\bz/n}(M)(S^{\e\prime})=\frac{\{(g,y)\in
G(S^{\e\prime})\times Y(S^{\e\prime})\colon
ng=-u(y)\}}{\{(-u(y),ny)\colon y\in Y(S^{\e\prime})\}}.
\end{equation}
The map $M\!\otimes^{\mathbf{L}}\! M^{*}\ra {\mathbb G}_{m}[1]$
induces a perfect pairing
$$
T_{\bz/n}(M)\times T_{\bz/n}(M^{*})\ra \mu_{\e n},
$$
where $\mu_{\e n}$ is the sheaf of $n$-th roots of unity. The above
pairing generalizes the classical Weil pairing of an abelian variety
$A$, which may be recovered by choosing $M=(0\ra A)$ and $n$ prime
to $p$ above. We will also need the following groups attached to
$M$:
$$
T_{p}(M)=\varprojlim T_{\bz/p^{\le m}}(M)
$$
(the {\it $p$-adic realization} of $M$), where the transition maps
are induced by the maps ${\s F}_{\,\bz/p^{\le m+1}}(M)\ra {\s
F}_{\,\bz/p^{\le m}}(M), [(g,y)]\mapsto [(p\e g,y)]$ (see (4)), and
$$
T(M)\{p\}=\varinjlim T_{\bz/p^{\le m}}(M)
$$
(the {\it $p$-divisible group attached to} $M$) with transition maps
induced by ${\s F}_{\,\bz/p^{\le m}}(M)\ra {\s F}_{\,\bz/p^{\le
m+1}}(M), [(g,y)]\mapsto [(g,p\e y)]$.

\smallskip

Now let $M$ be a 1-motive over $K$. For each prime $v$ of $K$, we
will write $M_{\e v}$ for the $K_{v}$-1-motive $M_{K_{v}}$. Further,
for each $i\geq -1$, ${\bh}^{\e i}(K_{v},M)$ will denote ${\bh}^{\e
i}(K_{v},M_{\e v})$. When $i=-1,1,2$, the group ${\bh}^{\e
i}(K_{v},M_{\e v})$ will be endowed with the discrete topology. For
$i=0$, ${\bh}^{\e i}(K_{v},M)$ will be endowed with the topology
defined in [10, p.99]. Define
$$
{\bh}^{\e -1}_{\,\wedge\be}(K_{v},M)=\text{Ker}\left[\,H^{\e
0}(K_{v},Y)\e\widehat{\,\phantom{.}}\ra H^{\e
0}(K_{v},G\e)\e\widehat{\,\phantom{.}}\,\e\right].
$$
Then there exists a surjective and continuous map of profinite
groups ${\bh}^{\e -1}(K_{v},M) \,\widehat{\,\phantom{.}}\ra{\bh}^{\e
-1}_{\,\wedge\be}(K_{v},M)$, and therefore an injection
$$
{\bh}^{\e -1}_{\,\wedge\be}(K_{v},M)^{D}\hookrightarrow (\e{\bh}
^{\e -1}(K_{v},M)
\,\widehat{\,\phantom{.}}\e)^{D}={\bh}^{\e -1}(K_{v},M)^{D}.
$$

\begin{theorem} There exists a continuous pairing
$$
{\bh}^{\e i}(K_{v},M)\times {\bh}^{\e 1-i}(K_{v},M^{\e *})
\ra{\bh}^{\e 2}(K_{v},{\mathbb G}_{m})=\bq/\e\bz
$$
which induces perfect dualities between the following profinite and
discrete groups, respectively.
\begin{enumerate}
\item[(a)]${\bh}^{\e -1}_{\,\wedge\be}(K_{v},M)$ and $\,{\bh}^{\e
2}(K_{v},M^{*})(p)$.
\item[(b)] ${\bh}^{\e
0}(K_{v},M)\e\widehat{\,\phantom{.}}$ and $\,{\bh}^{\e
1}(K_{v},M^{*})(p)$.
\end{enumerate}
For $i\neq -1,0,1,2$, the pairing is trivial.
\end{theorem}
\begin{proof} See [10, Theorem 2.3 and Lemma 2.1].
\end{proof}

\section{The Poitou-Tate exact sequence for $p$-primary finite and
flat group schemes}

Let $v$ be a prime of $K$ and let $N_{v}$ be a $p$-primary finite
and flat group scheme over $K_{v}$. The group $H^{\e i}(K_{v},N):=
H^{\e i}(K_{v},N_{v})$ (for $i=0,1$ or $2$) is canonically endowed
with a locally compact topology (see [15, comments following
III.6.5, p.341]). Both $H^{\e 0}(K_{v},N)$ and $H^{\e 2}(K_{v},N)$
are discrete (in fact, finite), but $H^{\e 1}(K_{v},N)$ need not be.
See [15, Table III.6.8, p.343]. We let $N_{v}^{d}={\s
Hom}(N_{v},{\bg}_{m})$ be the Cartier dual of $N_{v}$. Assume now
that $N_{v}$ extends to a finite and flat group scheme ${\s N}_{\e
v}$ over $\text{Spec}\,{\s O}_{v}$. Set $H^{\le i}\be\big(\e{\s
O}_{v},{\s N}\,\big)=H^{\le i}\be\big(\e{\s O}_{v},{\s
N}_{v}\,\big)$. By [13, p.293] or [15, beginning of \S III.7], for
each $i\geq 0$ the canonical map $H^{\le i}\be\big(\e{\s O}_{v},{\s
N}\,\big)\ra H^{\le i}(K_{v},N)$ embeds $H^{\le i}\be\big(\e{\s
O}_{v},{\s N}\,\big)$ as a compact and open subgroup of $H^{\le
i}(K_{v},N)$. We will identify $H^{\le i}\be\big(\e{\s O}_{v},{\s
N}\,\big)$ with its image in $H^{\le i}(K_{v},N)$ under this map.
Further, we will write $H^{\le i}_{v}({\s O}_{v},\s N)$ for the
cohomology group of $\s N$ with support on $\spec k(v)$ (see [15,
Proposition III.0.3, p.270]).

\begin{theorem} Let $i=0,1$ or $2$.
\begin{enumerate}
\item[(a)] There exists a perfect continuous pairing
$$
H^{\le i}(K_{v},N)\times H^{\le 2-i}(K_{v},N^{d})\ra
\bq_{p}/\e\bz_{\be p}.
$$
\item[(b)] In the pairing of (a), $H^{\le i}(\s O_{v},\s N\e)$ is
the exact annihilator of $H^{\le 2-i}(\s O_{v},\s N^{\e d})$.
\end{enumerate}
\end{theorem}
\begin{proof} For (a), see [15, Theorem III.6.10, p.344]. Statement
(b) is [op.cit., Corollary III.7.2, p.349].
\end{proof}

Now let $N$ be a $p$-primary finite and flat group scheme over $K$.
For any prime $v$ of $K$, we will write $N_{v}=N\times_{\spec K}
\spec K_{v}$ and $H^{\le i}(K_{v},N)=H^{\le i}(K_{v},N_{v})$. Recall
the set $\s F$ defined in Section 2. The elements of $\s F$ are
pairs $(U,\s N\e)$, where $U$ is a nonempty open affine subscheme of
$X$ such that $N$ extends to a finite and flat group scheme ${\s N}$
over $U$. If $(U,\s N\e)\in\s F$ and $v\in U$, we will write ${\s
N}_{v}= {\s N}\times_{U}\e\spec{\s O}_{v}$ and $H^{\e
i}\be\big(\e{\s O}_{v},{\s N}\e\big)=H^{\e i}\be\big(\e{\s
O}_{v},{\s N}_{v}\big)$. For $(U,\s N\e)\in\s F$, let $H^{\e
i}_{\textrm{c}}(U,{\s N}\e)$ be defined as in [15, comments
preceding Proposition III.0.4, pp.270-271] with $K_{v}$ replacing
the field of fractions of the henselization of the local ring at
$v$. Then [op.cit., Proposition III.0.4(a), p.271] remains valid,
i.e., there exists an exact sequence
\begin{equation}
\dots\ra H^{\e i}_{\e\text{c}}(U,\s N\e)\ra H^{\e i}(U,\s N\e)\ra
\displaystyle\bigoplus_{v\notin U}H^{\e i}(K_{v},N)\ra H^{\e
i+1}_{\e\text{c}}(U,\s N\e)\ra\dots.
\end{equation}
See [15, III.0.6(b), p.274]. The map $H^{\e i}(U,\s
N\e)\ra\bigoplus_{\e v\notin U} H^{\e i}(K_{v},N)$ appearing above
is the sum over $v\notin U$ of the maps $H^{\e i}(U,\s N\e)\ra H^{\e
i}(K_{v},N)$ induced by the composite morphism $\spec K_{v}\ra\spec
K\ra U$. Set
$$\begin{array}{rcl}
D^{\e i}(U,\s N\e)&=& \ng\text{Ker}\ng \left[H^{\e i}(U,\s
N\e)\ra\displaystyle\bigoplus_{v\notin U}
H^{\e i}(K_{v},N)\right]\\\\
&=&{\rm{Im}}\be\left[\,H^{\e i}_{\e\rm{c}}(U,\s N\e)\to H^{\e
i}(U,\s N\e)\,\right].
\end{array}
$$
Since $H^{\e 0}(U,\s N\e)=N(K)$ and $H^{\e 0}(K_{v},N)=N(K_{v})$ for
every $v\notin U$, the map $H^{\e 0}(U,\s N\e)\ra\bigoplus_{\e
v\notin U} H^{\e i}(K_{v},N)$ is injective, i.e., $D^{\e 0}(U,\s
N\e)=0$.

\smallskip

\begin{lemma} For any $(U,\s N\e)\in\s F$, the canonical map
$H^{\e 1}(U,\s N\e)\ra H^{\e 1}(K,N)$ is injective.
\end{lemma}
\begin{proof} The proof is similar to the proof of [15, Lemma
III.1.1, p.286].
\end{proof}

Using the above lemma, we will regard $D^{\e 1}(U,\s N\e)$ as a
subgroup of $H^{\e 1}(K,N)$ for any $(U,\s N\e)\in\s F$.

From now on, we will simplify our notations by writing $(V,\s N\e)$
for $(V,\s N\!\be\mid_{\e V}\be)$ when $(U,\s N\e)\in\s F$ and $V$
is an open subset of $U$.

\begin{lemma} Let $(U,\s N\e)\in\s F$ be arbitrary. Then there exists
a nonempty open subset $U_{0}$ of $U$ such that, for any nonempty
open subset $V\subset U_{0}$, both $D^{\e 1}(V,\s N)$ and $D^{\e
1}(V,\s N^{\e d})$ are finite.
\end{lemma}
\begin{proof} By a theorem of
M.Raynaud (see [17] or [3, Theorem 3.1.1, p. 110]), there exist a
nonempty open subset $U_{0}(\s N)\subset U$, abelian $U_{0}(\s
N\e)$-schemes $\s A$ and $\s B$ and an exact sequence $0\ra \s
N_{0}\ra \s A\ra \s B\ra 0$, where $\s N_{0}=\s N\!\mid_{\, U_{0}(\s
N)}$ and the first nontrivial map is a closed immersion. Let $V$ be
any nonempty open subset of $U_{0}(\s N\e)$. Then $0\ra \s
N_{0}\be\!\mid_{\e V}\ra \s A\be\!\mid_{\e V}\ra \s B\be\!\mid_{\e
V}\ra 0$ and $0\ra N_{K_{v}}\ra A_{K_{v}}\ra B_{K_{v}}\ra 0$, for
any prime $v$ of $K$, are also exact. Here $A$ and $B$ denote,
respectively, the generic fibers of $\s A$ and $\s B$. Using these
exact sequences and the fact that $\s B\e(V)=B(K)$, we obtain an
exact commutative diagram
\[
\xymatrix{B(K)\ar[d]\ar[r]& H^{\e 1}(V,\s N_{0}\e)\ar[d]\ar[r] &
H^{\e 1}(V,\s A)_{\e m}\ar[d]\\
\displaystyle\bigoplus_{v\notin V} B(K_{v})\ar[r]& \displaystyle
\bigoplus_{v\notin V}H^{\e 1}(K_{v},N)\ar[r] &
\displaystyle\bigoplus_{v\notin V} H^{\e 1}(K_{v},A),
\\
}
\]
where $m$ is any integer which annihilates $\s N_{0}$. Since the
image of $B(K)$ in $H^{\e 1}(V,\s N_{0})$ is finite by the
Mordell-Weil theorem, the finiteness of $D^{\e 1}(V,\s N_{0})$
follows from that of $D^{\e 1} (V,\s A)_{\e m}=\Sha^{\e 1}(K,A)_{\e
m}$, which is the main theorem of [13] (for the last equality, see
[15, Lemma II.5.5, p.246]). Now repeat the proof with $\s N^{\e d}$
in place of $\s N$ and take $U_{0}=U_{0}(\s N\e)\cap U_{0}(\s N^{\e
d})$.
\end{proof}

We now define, for $i=1$ or $2$,
$$
\Sha^{\e i}(K,N)=\text{Ker}\!\left[\e H^{\e i}(K,N)\ra\displaystyle
\prod_{\text{all $v$}}H^{\e i}(K_{v},N)\right].
$$

\begin{proposition} Let $(U,\s N\e)\in\s F$ be arbitrary and let
$U_{0}$ be as in the statement of the previous lemma. Then there
exists a nonempty open subset $U_{1}\subset U_{0}$ such that, for
any nonempty open subset $V$ of $\,U_{1}$, $\!\Sha^{\e 1}(K,N)\be
=\be D^{\e 1}(V,\s N)$. In particular, $\Sha^{\e 1}(K,N)$ is a
finite group.
\end{proposition}
\begin{proof} (Cf. [10, proofs of Lemma 4.7 and Theorem 4.8,
pp.114-115]). By definition, $\Sha^{\e 1}(K,N)\supset\bigcap_{\,\,
\emptyset\neq W\subset U_{0}} D^{\e 1}(W,\s N\e)$. Since each set
$D^{\e 1}(W,\s N\e)$ is finite, we may choose finitely many nonempty
open subsets $W_{1},W_{2},\dots,W_{r}$ of $U_{0}$ such that
$$
\Sha^{\e 1}(K,N)\supset\displaystyle\bigcap_{j=1}^{r}D^{\e 1}
(W_{j},\s N\e).
$$
Let $U_{1}=\bigcap_{\e j=1}^{\e r}W_{j}$ and let $V$ be any nonempty
open subset of $U_{1}$. By [15, Proposition III.0.4(c), p.271, and
Remark III.0.6(b), p.274], for any $j$ there exist natural maps
$H_{\text{c}}^{\e 1}(V,\s N\e)\overset{f_{j}}\longrightarrow
H_{\text{c}}^{\e 1}(W_{j},\s N\e)\overset{g_{j}}\longrightarrow
H^{\e 1}(K,N)$ such that $\text{Im}\e(\e g_{j}\e\circ\e f_{j}\e)=
D^{\e 1}\lbe(V,\s N)$ and $\text{Im}\e(\e g_{j}\e)=D^{\e 1}(W_{
j},\s N)$. It follows that $D^{\e 1}(V,\s N\e)\subset D^{\e
1}(W_{j},\s N\e)$ for every $j$ and we conclude that $D^{\e 1}(V,\s
N\e)\subset \Sha^{\e 1}(K,N)$. To prove the reverse inclusion, let
$\xi\in \Sha^{\e 1}(K,N)$. Then $\xi$ extends to $H^{\e 1}(W,\s
N\e)$ for some nonempty open subset $W$ of $U$, which we may assume
to be contained in $V$. Then $\xi\in D^{\e 1}(W,\s N\e)\subset D^{\e
1}(V,\s N\e)$ (by the same argument as above), and the proof is
complete.
\end{proof}

\begin{lemma} Let $(U,\s N\e)\leq (V,\s N\e)
\in\s F$. Then the natural map $H^{\e 2}(U,\s N\e)\ra H^{\e 2}(V,\s
N\e)$ induces a map $D^{\e 2} (U,\s N\e)\ra D^{\e 2}(V,\s N\e)$.
\end{lemma}
\begin{proof} For each $v$, the boundary map $H^{\e 2}(K_{v},N)\ra
H^{\e 3}_{v}(\s O_{v},\s N)$ appearing in the localization sequence
for the pair $\spec K_{v}\subset\spec\s O_{v}$ [15, Proposition
III.0.3(c), p.270] is an isomorphism [op.cit., comments preceding
Theorem III.7.1, p.349]. Thus the localization sequence for the pair
$V\subset U$ induces an exact sequence
$$
H^{\e 2}(U,\s N\e)\ra H^{\e 2}(V,\s N\e)\ra
\displaystyle\bigoplus_{v\in U\setminus V} H^{\e 2} (K_{v},N).
$$
It is not difficult to check that the second map in the above exact
sequence is the natural one, from which the lemma follows.
\end{proof}

\begin{proposition} There exists a canonical isomorphism
$$
\displaystyle\varinjlim_{(U,\,\s N\e)\in\s F}D^{2}(U,\s N\e)=
\Sha^{\e 2}(K,N).
$$
\end{proposition}
\begin{proof} For any $(U,\s N\e)\in\s F$, set
$$
{\s D}^{\e 2}(U,\s N\e)=\text{Im}[\e D^{\e 2}(U,\s N\e)\ra H^{\e 2}
(K,N)\e].
$$
Let $(V,\s N\e)\in\s F$ be such that $(U,\s N\e) \leq (V,\s N\e)$.
By Lemma 4.5, the map $D^{\e 2}(U,\s N\e)\ra H^{\e 2}(K,N)$ factors
through ${\s D}^{\e 2}(V,\s N\e)$, whence ${\s D}^{\e 2}(U,\s
N\e)\subset {\s D}^{\e 2}(V,\s N\e)$. Now the identification of
$D^{\e 2}(U,\s N\e)$ with $\text{Im}[\e H^{\e 2}_{\text{c}}(U,\s
N\e)\ra H^{\e 2}(U,\s N\e)\e]$ and the covariance of $H^{\e
2}_{\text{c}}(-,\s N\e)$ with respect to open immersions show that
${\s D}^{\e 2}(V,\s N\e\e)\subset {\s D}^{\e 2}(U,\s N\e)$. We
conclude that ${\s D}^{\e 2}(V,\s N\e)={\s D}^{\e 2}(U,\s N\e)$ for
all $(V,\s N)$ as above and necessarily ${\s D}^{\e 2}(U,\s
N\e)=\Sha^{\e 2}(K,N)$ for any $(U,\s N\e)\in\s F$. Thus we have a
surjection
$$
\displaystyle\varinjlim_{(U,\,\s N\e)\in\s F}D^{2}(U,\s N\e)\ra
\Sha^{\e 2}(K,N).
$$
By Lemma 2.3 this is an injection as well, which completes the
proof.
\end{proof}

We would like to establish and analogue of the previous proposition
for $\Sha^{\e 1}(K,N)$. However (cf. Lemma 4.5), the natural map
$H^{\e 1}(U,\s N\e)\ra H^{\e 1}(V,\s N\e\e)$ need not map $D^{\e
1}(U,\s N\e)$ into $D^{\e 1}(V,\s N\e\e)$ for $(U,\s N\e)\leq(V,\s
N\e)\in\s F$. The reason for this is that a class $\xi\in D^{\e
1}(U,\s N\e)$ need not map to zero in $H^{\e 1}(K_{v},N)$ for primes
$v\in U\setminus V$. Following [10], we will circumvent this
difficulty in the next proposition by showing that the groups $D^{\e
1}(U,\s N\e)$ ``eventually become constant with value
$\!\be\be\Sha^{\e 1}(K,N)$'', by which we mean that there exists an
element $(U_{1},\s N\e)\in \s F$ such that, for every $(V,\s
N\e)\in\s F$ with $(U_{1},\s N\e)\leq (V,\s N\e)$, $D^{\e 1}(V,\s
N\e)$ can be identified with $\!\!\Sha^{\e 1}(K,N)$.

\begin{proposition} Let $(U,\s N\e)\in\s F$ be arbitrary and let
$U_{0}\subset U$ be as in the statement of Lemma 4.3. Then there
exists a nonempty open subset $U_{1}\subset U_{0}$ such that, for
any nonempty open subset $V$ of $U_{1}$, $\!\Sha^{\e 1}(K,N)\be =\be
D^{\e 1}(V,\s N\e)$. In particular, $\Sha^{\e 1}(K,N)$ is a finite
group.
\end{proposition}
\begin{proof} (Cf. [10, proofs of Lemma 4.7 and Theorem 4.8,
pp.114-115]). By definition, $\Sha^{\e 1}(K,N)\supset\bigcap_{\,\,
\emptyset\neq W\subset U_{0}} D^{\e 1}(W,\s N\e)$. Since each set
$D^{\e 1}(W,\s N\e)$ is finite (see Lemma 4.3), we may choose
finitely many nonempty open subsets $W_{1},W_{2},\dots,W_{r}$ of
$U_{0}$ such that
$$
\Sha^{\e 1}(K,N)\supset\displaystyle\bigcap_{j=1}^{r}D^{\e 1}
(W_{\be j},\s N\e).
$$
Let $U_{1}=\bigcap_{\e j=1}^{\e r}W_{\be j}$ and let $V$ be any
nonempty open subset of $U_{1}$. By [15, Proposition III.0.4(c),
p.271, and Remark III.0.6(b), p.274], for any $j$ there exist
natural maps $H_{\text{c}}^{\e 1}(V,\s
N\e)\overset{f_{j}}\longrightarrow H_{\text{c}}^{\e 1}(W_{j},\s
N\e)\overset{g_{j}}\longrightarrow H^{\e 1}(K,N)$ such that
$\text{Im}\e(\e g_{j}\e\circ\e f_{j}\e)= D^{\e 1}\be(V,\s N\e)$ and
$\text{Im}\e(\e g_{j}\e)=D^{\e 1}(W_{ j},\s N\e)$. It follows that
$D^{\e 1}(V,\s N\e)\subset D^{\e 1}(W_{j},\s N\e)$ for every $j$ and
we conclude that $D^{\e 1}(V,\s N\e)\subset \Sha^{\e 1}(K,N)$. To
prove the reverse inclusion, let $\xi\in \Sha^{\e 1}(K,N)$. Then
$\xi$ extends to $H^{\e 1}(W,\s N\e)$ for some nonempty open subset
$W$ of $U$, which we may assume to be contained in $V$. Then $\xi\in
D^{\e 1}(W,\s N\e)\subset D^{\e 1}(V,\s N\e)$ (by the same argument
as above), and the proof is complete.
\end{proof}

Let $(U,\s N)$ be arbitrary and let $U_{1}$ be as in the statement
of the previous proposition. Set
$$
\s F_{1}=\{(V,\s N\e)\in\s F\colon (U_{1}, \s N)\leq (V,\s N\e)\}.
$$
Then, if $(V,\s N\e)\leq (W,\s N\e)\in\s F_{1}$, there exist natural
maps $D^{\e 1}(W,\s N\e)\ra D^{\e 1} (V,\s N\e)$ (the identity map;
see Proposition 4.7) and $D^{\e 2}\!\left(V,{\s N}^{\e
d}\e\right)\ra D^{\e 2}\!\left(W,{\s N}^{\e d}\e\right)$ (see Lemma
4.5). The respective limits are
\begin{equation}
\displaystyle\varprojlim_{\,(V\be,\,\e\s N\e) \,\in\,\s F_{1}}\be
D^{\e 1}(V,\s N\e)= \Sha^{\e 1}(K,N)
\end{equation}
and
\begin{equation}
\displaystyle\varinjlim_{\,(V\be,\,\e\s N\e) \,\in\,\s F_{1}}\be
D^{\e 2}\!\be\left(V, (\s N)^{d}\e\right)=\Sha^{\e
2}\be\!\left(K,N^{d}\e \right)
\end{equation}
(see Proposition 4.6).

\begin{lemma} Let $(U,\s N\e)\in\s F$ be arbitrary and let
$U_{1}\subset U$ be as in the statement of Proposition 4.7. Then,
for any nonempty open subset $V\subset U_{1}$ and $i=1$ or 2, there
exists a perfect pairing of finite groups
$$
D^{\e i}(V,\s N\e)\times D^{\e 3-i}(V,{\s N}^{\e
d}\e)\ra\bq_{p}/\e\bz_{\be p}\,.
$$
\end{lemma}
\begin{proof} By [15, Theorem III.8.2, p.361], for any $i$ there
exists a perfect pairing
$$
[-,-]\colon H^{\e i}(V,\s N\e)\times H^{\e 3-i}_{\e\text{c}}(V,{\s
N}^{\e d}\e)\ra\bq_{p}/\e\bz_{\be p}
$$
between the torsion discrete group $H^{\e i}(V,\s N\e)$ and the
profinite group $H^{\e 3-i}_{\e\text{c}}(V,{\s N}^{\e d}\e)$. The
above pairing induces the middle vertical map in the following
natural commutative diagram:
\[
\xymatrix{0\ar[r]& D^{\e i}(V,\s N\e)\ar[r] & {H}^{\e i}(V,\s
N\e)\ar[d]\ar[r] & \displaystyle\bigoplus_{v\notin V} {H}^{\e
i}(K_{v},N)\ar[d]\\
&& {H}^{\e 3-i}_{\e\text{c}} (V,\s N^{\e d}\e)^{D}\ar[r] &
\displaystyle\bigoplus_{v\notin V} {H}^{\e 2-i}(K_{v},N^{\e
d})^{D}.}
\]
It follows that there exists a well-defined pairing
$$
(-,-)\colon D^{\e i}(V,\s N\e)\times D^{\e 3-i}(V,{\s N}^{\e
d}\e)\ra\bq_{p}/\e\bz_{\be p}
$$
given by $(a,a^{\e\prime})=[\e a,b^{\e\prime}\e]$, where $a\in D^{\e
i}(V,\s N\e)\subset H^{\e i}(V,\s N\e)$ and $b^{\e\prime}$ is a
preimage of $a^{\e\prime}$ in $H^{\e 3-i}_{\e\text{c}}(V,{\s N}^{\e
d}\e)$. The above pairing is non-degenerate on the left. Similarly,
interchanging $i$ and $3-i$ and $\s N$ and $\s N^{\e d}$ above, we
obtain a pairing
$$
D^{\e 3-i}(V,\s N^{\e d}\e)\times D^{\e i}(V,{\s
N}\e)\ra\bq_{p}/\e\bz_{\be p}
$$
with trivial left kernel. Setting $i=1$ above and using the
finiteness of $D^{\e 1}(V,{\s N}\e)$ (see Lemma 4.3 and recall that
$U_{1}\subset U_{0}$), we obtain the case $i=1$ of the lemma. To
obtain the case $i=2$ of the lemma, one argues similarly,
interchanging the roles of $\s N$ and $\s N^{\e d}$ and using the
finiteness of $D^{\e 1}(V,{\s N}^{\e d}\e)$ (see Lemma 4.3).
\end{proof}

\begin{theorem} There exists a perfect pairing of finite groups
$$
\Sha^{\e 1}(K,N)\times\!\! \Sha^{\e 2}\!\be\left(K,N^{\e d}\e\right)
\ra\bq_{\le p}/\e\bz_{\be p}.
$$
\end{theorem}
\begin{proof} This follows at once from Lemma 4.8, using (6), (7)
and the finiteness of $\!\Sha^{\e 1}(K,N)$ (see Proposition 4.4).
\end{proof}

We will need the following lemma.

\begin{lemma} Let $U_{1}$ be as in the statement of Proposition 4.7.
Then there exists a canonical isomorphism
$$
H^{\e 2}(K,N^{d})^{D}=\displaystyle\varprojlim_{V\subset U_{1}}H^{\e
2}(V,\s N^{\e d}\e)^{D}.
$$
\end{lemma}
\begin{proof} Let $V$ be a nonempty open subset of $U_{1}$. Since
$D^{\e 2}(V,\s N^{\e d})$ is finite by Lemma 4.8, $H^{\e 2}(V,\s
N^{\e d}\e)$ is finite as well since it fits into an exact sequence
$$
0\ra D^{\e 2}(V,\s N^{\e d})\ra H^{\e 2}(V,\s N^{\e
d}\e)\ra\displaystyle\bigoplus_{v\notin V}H^{\e 2}(K_{v},N^{d}).
$$
Therefore $\varprojlim_{\,V\subset U_{1}}H^{\e 2}(V,\s N^{\e
d}\e)^{D}$ is canonically isomorphic to the dual of
$\varinjlim_{\,V\subset U_{1}}H^{\e 2}(V,\s N^{\e d}\e)$. Now Lemma
2.3 completes the proof.
\end{proof}

\medskip

For $(U,\s N\e)\in\s F$ and $0\leq i\leq 2$, define
$$
P^{\, i}(U,\s N\e)=\bigoplus_{v\notin U}H^{\le i}(K_{v},N)\times
\prod_{v\in U}H^{\le i}\be\big(\e{\s O}_{v},{\s
N}\,\big)\subset\displaystyle\prod_{\text{all $v$}}H^{\le
i}(K_{v},N)
$$
with the product topology. It is a locally compact group. Now, for
every $v$, $H^{\le 0}\be\big(\e{\s O}_{v},{\s N}\,\big)=H^{\le
0}(K_{v},N)$ and $H^{\e 2}\be\big(\e{\s O}_{v},{\s N}\e\big)=0$ [15,
beginning of III.7, p. 348], whence
$$
P^{\, 0}(U,\s N\e)=\prod_{\text{all $v$}}H^{\e 0}(K_{v},N)
$$
and
$$
P^{\, 2}(U,\s N\e)=\bigoplus_{v\notin U}H^{\e 2}(K_{v},N).
$$
Note that $P^{\, 0}(U,\s N\e)$ is compact and $P^{\, 2}(U,\s N\e)$
is finite. Further, if $(U,\s N\e)\leq (V,\s N\e)\in\s F$, then
$P^{\, i}(U,\s N\e)\subset P^{\, i}(V,\s N\e\e)$. Define
$$\begin{array}{rcl}
P^{\, i}(K,N)&=&\displaystyle\varinjlim_{(U,\e\s N\e)\in\s F}P^{\,
i}(U,\s N\e)\\\\
&=&\displaystyle\bigcup_{(U,\e\s N\e)\in\s F}P^{\, i}(U,\s
N\e)\subset\displaystyle\prod_{\text{all $v$}}H^{\le i}(K_{v},N),
\end{array}
$$
where the transition maps in the direct limit are the inclusion
maps. Thus $P^{\, i}(K,N)$ is the restricted topological product
over $v$ of the groups $H^{\le i}(K_{v},N)$ with respect to the
subgroups $H^{\le i}\be\big(\e{\s O}_{v},{\s N}\,\big)$. Now equip
$H^{\e i}(U,\s N\e)$ with the discrete topology. There exists a
natural map $H^{\e i}(U,\s N\e)\ra\prod_{\e v\in U}H^{\le
i}\be\big(\e{\s O}_{v},{\s N}\,\big)$, namely the product over $v\in
U$ of the maps $H^{\e i}(U,\s N\e)\ra H^{\le i}\be\big(\e{\s
O}_{v},{\s N}\,\big)$ induced by the canonical morphisms
$\spec{\mathcal O}_{v}\ra U$. The product of the above map with the
map $H^{\e i}(U,\s N\e)\ra \prod_{\e v\notin U}H^{\e i}(K_{v},N)$
introduced previously is a map
$$
\beta_{\e i}(U,\s N\e)\colon H^{\e i}(U,\s N\e)\ra P^{\, i}(U,\s
N\e).
$$
If $(U,\s N\e)\leq (V,\s N\e\e)\in\s F$, then there exists a
canonical commutative diagram
\[
\xymatrix{H^{\e i}(U,\s N\e)\ar[d]\ar[r]^{\beta_{\e i}(U,\s N\e)}&
P^{\e i}(U,\s
N\e)\ar@{^{(}->}[d]\\
H^{\e i}(V,{\s N}\e)\,\,\ar[r]^{\,\beta_{\e i}(V\be,\e\s N\e)\,}&
\,\,P^{\e i}(V,{\s N}\e),}
\]
where the left-hand vertical map is induced by the inclusion
$V\subset U$. Consequently, by Lemma 2.3, the direct limit of the
maps $\beta_{\e i}(U,\s N\e)$ is a map
$$
\beta_{\e i}(K,N)\colon H^{\e i}(K,N)\ra P^{\, i}(K,N)
$$
whose kernel is $\Sha^{\e i}(K,N)$. Note that $\beta_{\e 0}(K,N)$ is
injective since it coincides with the canonical map
$N(K)\ra\prod_{\,\text{all $v$}} N(K_{v})$.

Now Theorem 4.1 shows that $P^{\, i}(K,N)$ is the algebraic and
topological dual of $P^{\,2-i}(K,N^{\e d})$. Indeed, there exists an
isomorphism of topological groups
$$
\phi^{\e i}_{(K,N)}\colon P^{\, i}(K,N)\ra P^{\,2-i}(K,N^{\e d})^{D}
$$
defined as follows. If $(\xi_{v})\in P^{\, i}(K,N)$ and
$(\zeta_{v})\in P^{\, 2-i}(K,N^{\e d})$ is arbitrary, then
$$
\phi^{\e i}_{(K,N)}(\xi_{v})(\zeta_{\e
v})=\displaystyle\sum_{\text{all $v$}} \,\,(\xi_{v},\zeta_{\e
v})_{v}\,\in\,\bq_{p}/\e\bz_{\be p}\,,
$$
where, for each $v$, $(-,-)_{v}$ is the pairing of Theorem 4.1(a)
(that the sum is finite follows from the definition of $P^{\,
i}(K,N)$ and the fact that, for each $v$, $H^{\le i}\be\big(\e{\s
O}_{v},{\s N}\,\big)$ and $H^{\le 2-i}\be\big(\e{\s O}_{v},{\s
N}^{\e d}\,\big)$ annihilate each other under the pairing
$(-,-)_{v}$). Now, for each $(U,\s N)\in\s F$, let
$$
\phi^{\e i}_{(U,\e\s N\e)}\colon P^{\e i}(U,\s N\e)\ra P^{\e 2-
i}(U,\s N^{\e d}\e)^{D}
$$
be the composition of the restriction of $\phi^{\e i}_{(K,N)}$ to
$P^{\e i}(U,\s N\e)$ and the canonical map $P^{\,2-i}(K,N^{\e
d})^{D}\ra P^{\e 2- i}(U,\s N^{\e d}\e)^{D}$. Note that, if
$(\xi_{v})\in P^{\, i}(U,\s N\e)$ and $(\zeta_{v})\in P^{\,
2-i}(U,\s N^{\e d})$ is arbitrary, then
$$
\phi^{\e i}_{(U,\e\s N\e)}(\xi_{v})(\zeta_{\e
v})=\displaystyle\sum_{v\notin U} \,\,(\xi_{v},\zeta_{v})_{v}.
$$
We now let
$$
\gamma_{i}(U,\s N\e)\colon P^{\e i}(U,\s N\e)\ra H^{\e 2-i}(U,\s
N^{\e d}\e)^{D}
$$
be the composite $\beta_{\e 2-i}(U,{\s N}^{\e d}\e)^{D}\circ\phi^{\e
i}_{\e(U,\,\s N\e)}$. Further, let
$$
\psi^{\e i}_{\e (U,\e\s N\e)}\colon\displaystyle\bigoplus_{v\notin
U}H^{\le i}(K_{v},N)\ra H^{\e 2-i}(U,\s N^{\e d}\e)^{D}
$$
be given by
$$
\psi^{\e i}_{\e (U,\e\s N\e)}((\xi_{v}))(\zeta)=\sum_{v\notin
U}\,(\xi_{v},\zeta\!\be\mid_{\lbe K_{\lbe v}}\be)_{v}
$$
where, for each $v\notin U$, $\zeta\!\be\mid_{\lbe K_{\lbe v}}$
denotes the image of $\zeta\in H^{\e 2-i}(U,\s N^{\e d})$ in $H^{\e
2-i}(K_{v}, N^{d})$ under the map
$$
H^{\e 2-i}(U,\s N^{\e d})\ra H^{\e 2-i}(K, N^{d})\ra H^{\e
2-i}(K_{v}, N^{d})
$$
induced by the composite morphism $\spec K_{v}\ra\spec K\ra U$. Then
the following diagrams commute
\begin{equation}
\xymatrix{P^{\e i}(U,\s N)\ar[r]^(0.45){\gamma_{\e i}(U,\e\s N\e)}&
H^{\e
2-i}(U,\s N^{d})^{D}\\
\displaystyle\bigoplus_{v\notin U}H^{\le
i}(K_{v},N)\ar@{^{(}->}[u]\ar[ur]_(0.6){\psi^{\e i}_{(U,\e\s N)}}&&}
\end{equation}
and
\begin{equation}
\xymatrix{P^{\e i}(U,\s N)\ar[d]\ar[r]^(0.45){\gamma_{\e i}(U,\e\s
N\e)}& H^{\e
2-i}(U,\s N^{d})^{D}\\
\displaystyle\bigoplus_{v\notin U}H^{\le
i}(K_{v},N)\,,\ar[ur]_(0.6){\psi^{\e i}_{(U,\e\s N)}}&&}
\end{equation}
where the vertical map in (9) is the canonical projection. Now let
$$
\gamma_{i}(K,N)\colon P^{\, i}(K,N)\ra H^{\e 2-i}(K,N^{\e d})^{D}
$$
be the composite $\beta_{2-i}(K,N^{\e d})^{D}\be\circ\e \phi^{\e
i}_{(K,N)}$. Then, for any $(U,\s N\e)\in\s F$, $\gamma_{\e
i}(U,\e\s N\e)$ factors as
\begin{equation}
P^{\e i}(U,\s N)\hookrightarrow P^{\e
i}(K,N)\overset{\gamma_{i}(K,N)}\longrightarrow H^{\e
2-i}(K,N^{d})^{D}\ra H^{\e 2-i}(U,\s N^{\e d})^{D},
\end{equation}
where the last map is the dual of the canonical map $H^{\e 2-i}(U,\s
N^{\e d})\ra H^{\e 2-i}(K,N^{\e d})$ induced by the morphism $\spec
K\ra U$.

\begin{proposition} For $i=0,1$ or $2$, there exists a natural
isomorphism $\krn(\gamma_{\e i}(K,N))=\img(\e\beta_{\e i}(K,N))$.
Further, $\gamma_{\e 2}(K,N)$ is surjective.
\end{proposition}
\begin{proof} We show first that $\krn(\gamma_{\e i}(K,N\e))
\subset\img(\beta_{\e i}(K,N\e))$. Let $(\xi_{v})\in\krn(\gamma_{\e
i}(K,N\e))$. Then $(\xi_{v})\in P^{\e i}(U,\s N\e)$ for some
nonempty open affine subset of $X$, which we may assume to be
contained in the set $U_{1}$ introduced in the proof of Proposition
4.7. The element $(\xi_{v})_{v\notin U}\in\bigoplus_{\, v\notin
U}H^{\le i}(K_{v},N)$ is in the kernel of $\psi^{\e i}_{(U,\e\s
N\e)}$ by the commutativity of (8). Therefore $(\xi_{v})_{v\notin
U}$ is in the kernel of the composite map
\begin{equation}
\displaystyle\bigoplus_{v\notin U} H^{\e i}(K_{v},N)\overset{
\psi^{\e i}_{(U,\e\s N)}}\longrightarrow H^{\e 2-i}(U,\s N^{\e
d}\e)^{D}\overset{\sim}\longrightarrow H^{\e i+1}_{\e\text{c}}(U,\s
N\e),
\end{equation}
where the last isomorphism is induced by the perfect pairing
\begin{equation}
[-,-]\colon H^{\e 2-i}(U,\s N^{\e d}\e)\times H^{\e
i+1}_{\e\text{c}}(U,{\s N}\e)\ra\bq_{p}/\e\bz_{\be p}
\end{equation}
between the discrete torsion group $H^{\e 2-i}(U,\s N^{\e d}\e)$ and
the profinite group $H^{\e i+1}_{\e\text{c}}(U,{\s N}\e)$ (see [15,
Theorem III.8.2, p.361]). Consequently, by the exactness of (5),
there exists an element $\e\xi_{_{U}}\be\in H^{\e i}(U,\s N\e)$ such
that $\xi_{_{U}}\be\be\!\!\mid_{K_{v}}=\xi_{v}$ for all $v\notin U$.
The assignment $(\xi_{v})_{v\notin U}\mapsto\xi_{_{U}}$ is
functorial in $U$, i.e., $\xi_{_{\le V}}=\xi_{_{\le U}}\!\!\mid_{V}$
if $(U,\s N\e)\leq (V,\s N\e)\in\s F$, where $\xi_{_{\le
U}}\!\!\mid_{V}$ denotes the image of $\xi_{_{\le U}}$ under the map
$H^{\e i}(U,\s N\e)\ra H^{\e i}(V,\s N\e)$ induced by the inclusion
$V\subset U$. For any $U$ as above, let $\xi=\xi_{_{\le
U}}\!\!\mid_{K}$ be the image of $\xi_{_{U}}$ under the map $H^{\e
i}(U,\s N\e)\ra H^{\e i}(K,N)$ induced by the morphism $\spec K\ra
U$. Then $\xi$ is a well-defined element of $H^{\e i}(K,N)$ whose
image under $\beta_{\e i}(K,N\e)$ is $(\xi_{v})$. This shows that
$\krn(\gamma_{\e i}(K,N\e))\subset\img(\beta_{\e i}(K,N\e))$.

Next we will show that $\gamma_{\e i}(K,N\e)\circ \beta_{\e
i}(K,N\e)=0$, which will show that $\img(\beta_{\e
i}(K,N\e))\subset\krn(\gamma_{\e i}(K,N\e))$. Let $\xi\in H^{\e
i}(K,N)$. By Lemma 2.3, there exists a pair $(U,\s N\e)\in\s F$ and
an element $\xi_{_{U}}\in H^{\e i}(U,\s N\e)$ such that
$\xi=\xi_{_{\le U}}\be\!\!\mid_{K}$. For any $(V,\s N\e\e)\in\s F$
such that $(U,\s N\e)\leq (V,\s N\e)$, set $\xi_{_{\le
V}}=\xi_{_{\le U}}\!\!\mid_{V}$. The image of $\xi_{_{\le V}}$ under
the map $H^{\e i}(V,\s N\e)\ra\bigoplus_{v\notin V} H^{\e
i}(K_{v},N)$ maps to zero under the composite map (11) (with $U$
replaced by $V$ there). Consequently, it maps to zero under the map
$\psi^{\e i}_{(V,\e\s N\e)}$. Now the commutativity of (9) (with $U$
replaced by $V$ there) shows that $\gamma_{i}\big(V\be,\s
N\e\big)\big(\beta_{i}(V,\s N\e\big)\big(\xi_{_{V}}\big)\big)=0$. It
now follows from (10) (with $U$ replaced by $V$ there) that
$\gamma_{i}(K,N)(\beta_{i}(K,N)(\xi))$ is in the kernel of the
canonical map $H^{\e 2-i}(K,N^{\e d})^{D}\ra H^{\e 2-i}(V,\s N^{\e
d})^{D}$. This holds for any $V$ as above. Since the canonical map
$H^{\e 2-i}(K,N^{\e d})^{D}\ra\varprojlim_{V}H^{\e 2-i}(V,\s N^{\e
d})^{D}$ is an isomorphism (see Lemmas 2.4 and 4.10 and note that
this statement is trivially true if $i=2$), we conclude that
$\gamma_{i}(K,N)(\beta_{i}(K,N)(\xi))=0$, as desired.

It remains to show that $\gamma_{\e 2}(K,N)$ is surjective. In fact,
we will show that $\gamma_{\e 2}(U,\s N\e)$ is surjective for any
$(U,\s N\e)\in\s F$. By (10) and the bijectivity of the canonical
map $H^{\e 0}(K,N^{d})^{D}\ra H^{\e 0}(U,\s N^{\e d})^{D}$, this
will complete the proof. Since $``\,\text{Coker}\e(f)^{D}=
\text{Ker}\left(f^{D}\right)$'' if $f$ is a map between finite
groups, we have a canonical isomorphism
$$
{\rm{Coker}}\e(\gamma_{\e 2}(U,\s N\e))^{D}=
\text{Ker}\ng\!\left[\,H^{\e 0}(U,{\s N}^{\e
d}\e)\overset{\gamma_{\e 2}(U,\,\s N\e)^{D}}\longrightarrow P^{\,
2}(U,\s N\e)^{D}\right].
$$
The map $\gamma_{\e 2}(U,\,\s N\e)^{D}$ may be identified with the
natural map $N^{ d}(K)\ra\bigoplus_{\e v\notin U}N^{ d}(K_{v})$,
which is clearly injective. Thus $\gamma_{\e 2}(U,\,\s N\e)$ is
indeed surjective.
\end{proof}

Now, for $i=0$ or $1$, we define a map
$$
\delta_{\e i}(K,N)\colon H^{\e 2-i}\big(K,N^{d}\e\big)^{D}\ra H^{\e
i+1}(K,N)
$$
as the composite
$$
H^{\e 2-i}\big(K,N^{d}\e\big)^{D}\twoheadrightarrow\Sha^{\e
2-i}(K,N^{d})^{D}\simeq \Sha^{\e i+1}(K,N)\hookrightarrow H^{\e
i+1}(K,N),
$$
where the isomorphism comes from Theorem 4.9 (applied to $N$ and
$N^{d}$). Clearly
$$\begin{array}{rcl}
\krn(\delta_{\e i}(K,N))&=&\img\!\left[\e P^{\,2-i}(K,N^{\e
d})^{D}\ra H^{\e
2-i}\big(K,N^{d}\e\big)^{D}\right]\\
&=&\img\!\left[\e P^{\, i}(K,N)\ra H^{\e
2-i}\big(K,N^{d}\e\big)^{D}\right]\\
&=&\img(\gamma_{i}(K,N)).
\end{array}
$$
Further, $\img(\delta_{\e i}(K,N))=\Sha^{\e
i+1}(K,N)=\krn(\e\beta_{\e i+1}(K,N))$. These facts, together with
Proposition 4.11, yield the following {\it Poitou-Tate exact
sequence in flat cohomology}.

\begin{theorem} There exists an exact sequence of locally compact
groups and continuous homomorphisms
\[
\xymatrix{& 0\ar[r]& H^{\e 0}(K,N)\ar[r]^{\beta_{\e 0}} & P^{\e
0}(K,N)\ar[r]^{\gamma_{0}}& H^{\e 2}(K,N^{\e d})^{\e D}\ar[d]
^{\delta_{\e 0}}\\
&& \ar[d]^{\delta_{\e 1}} H^{\e 1}(K,N^{\e d})^{\e D} &
\ar[l]_{\gamma_{\e 1}}P^{\e 1}(K,N)
&\ar[l]_{\beta_{\e 1}} H^{\e 1}(K,N) \\
&& H^{\e 2}(K,N)\ar[r]^{\beta_{\e 2}}& P^{\e
2}(K,N)\ar[r]^{\gamma_{\e 2}}& H^{\e 0}(K,N^{\e d})^{\e D}\ar[r]&0.}
\]\qed
\end{theorem}

\section{1-Motives over open affine subschemes of $X$}

In this section all groups will be endowed with the discrete
topology, with the exception of the groups ${\bh}^{\e 0}(K_{v},M)$,
which will be endowed with the topology defined in [10, p.99]. Thus
${\bh}^{\e 0}(K_{v},M)$ contains a closed subgroup of finite index
which is a (possibly non-Hausdorff) quotient of $G(K_{v})$. The
theory of Lie groups over a local field shows that $G(K_{v})$ is
locally compact, compactly generated and completely disconnected.
Therefore, by [12, Theorem II.9.8, p.90], $G(K_{v})$ is
topologically isomorphic to a product $\bz^{\e b}\times C$, where
$b$ is a non-negative integer and $C$ is a compact and completely
disconnected, i.e., profinite, abelian group. We conclude that, if
${\bh}^{\e 0}(K_{v},M)^{\e\prime}$ denotes the quotient of
${\bh}^{\e 0}(K_{v},M)$ by the closure of $\{0\}$ and $n$ is any
integer, then ${\bh}^{\e 0}(K_{v},M)^{\e\prime}/n$ is a profinite
group.

Let $U$ be any nonempty open affine subscheme of $X$. For any
cohomologically bounded complex ${\s F}^{\e\bullet}$ of fppf sheaves
on $U$, there exist cohomology groups with compact support
${\bh}^{\e i}_{\e\text{c}}(U,{\s F}^{\e\bullet}\e)$ which may be
defined as in [15, comments preceding Proposition III.0.4, p.271].
There exists an exact sequence
$$
\dots\ra\! {\bh}^{\e i}_{\e\text{c}}(U,{\s F}^{\e\bullet}\e)\!\ra\!
{\bh}^{\e i}(U,{\s F}^{\e\bullet}\e)\!\ra\!
\displaystyle\bigoplus_{v\notin U}{\bh}^{\e i}(K_{v},{\s
F}^{\e\bullet}\e)\!\ra\! {\bh}^{\e i+1}_{\e\text{c}}(U,{\s
F}^{\e\bullet}\e)\!\ra\dots,
$$
where we have abused notation in the third term by writing ${\s
F}^{\e\bullet}$ for the pullback of ${\s F}^{\e\bullet}$ under the
composite map $\spec K_{v}\ra\spec K\ra U$.

For any pair of cohomologically bounded complexes ${\s
F}^{\e\bullet}$ and ${\s G}^{\e\bullet}$ as above, there exists a
cup-product pairing
$$
{\bh}^{\e i}(U,{\s F}^{\e\bullet})\times{\bh}^{\e
j}_{\e\text{c}}(U,{\s G}^{\e\bullet})\ra {\bh}^{\e
i+j}_{\e\text{c}}(U,{\s F}^{\e\bullet}\!\otimes^{\,\mathbf{L}} {\s
G}^{\e\bullet}).
$$

Now let $\s M$ be a 1-motive over $U$. Set
$$
H^{\e i}(U,T_{p}(\s M))=\varprojlim_{m} H^{\e i}(U,T_{\,\bz/\e
p^{\le m}}(\s M))
$$
and
$$
H^{\e i}_{\e{\rm{c}}}(U,T_{p}(\s M))=\varprojlim_{m} H^{\e
i}_{\e{\rm{c}}}(U,T_{\,\bz/\e p^{\le m}}(\s M)).
$$

\begin{lemma} Let $i=0,1$ or $\e 2$.
\begin{enumerate}
\item[(a)] There exists a pairing
$$
H^{\e i+1}\be\left(U,T_{p}(\s M)\right)(p)\times (H^{\e 2-
i}_{\e{\rm{c}}}(U,T(\s M^{*})\{p\})/\e
p{\rm{-div}})\ra\bq_{p}/\e\bz_{\be p}
$$
which is non-degenerate on the left.
\item[(b)] There exists a pairing
$$
(H^{\e i}\be\left(U,T(\s M)\{p\}\right)/\e p{\rm{-div}})\times H^{\e
3- i}_{\e{\rm{c}}}(U,T_{p}(\s M^{*}))(p)\ra\bq_{p}/\e\bz_{\be p}
$$
which is non-degenerate on the right.
\end{enumerate}
\end{lemma}
\begin{proof} By [15, Theorem III.8.2, p.361],
for every $r\geq 0$ and any $m\geq 1$ there exists a perfect pairing
$$
H^{\e r}\!\left(U,T_{\,\bz/\e p^{\le m}}(\s M)\right)\times H^{\e
3-r}_{\text{c}}\be\left(U,T_{\,\bz/\e p^{\le m}}(\s
M^{*})\right)\ra\bq_{p}/\e\bz_{\be p}
$$
between the discrete torsion group $H^{\e r}\!\left(U,T_{\,\bz/\e
p^{\le m}}(\s M)\right)$ and the profinite group $H^{\e
3-r}_{\text{c}}\be\left(U,T_{\,\bz/\e p^{\le m}}(\s M^{*})\right)$.
Setting $r=i+1$ and $r=i$ above, we obtain perfect pairings
$$
H^{\e i+1}\be\left(U,T_{p}(\s M)\right)\times H^{\e 2-
i}_{\e{\rm{c}}}(U,T(\s M^{*})\{p\})\ra\bq_{p}/\e\bz_{\be p}
$$
and
$$
H^{\e i}\be\left(U,T(\s M)\{p\}\right)\times H^{\e 3-
i}_{\e{\rm{c}}}(U,T_{p}(\s M^{*}))\ra\bq_{p}/\e\bz_{\be p}.
$$
The lemma now follows easily.
\end{proof}

\smallskip

For each $i$ such that $-1\leq i\leq 3$, there exists a canonical
pairing
\begin{equation}
\langle-,-\rangle\,\colon{\bh}^{\e i}(U,{\s M})\times{\bh}^{\,
2-i}_{\e\text{c}}(U,{\s M}^{*})\ra \bq/\e\bz.
\end{equation}
See [10, p.108]. The above pairing induces a pairing
\smallskip
\begin{equation}
{\bh}^{\e i}(U,{\s M})(p)/\pdiv\times{\bh}^{\,
2-i}_{\e\text{c}}(U,{\s M}^{*})(p)/\pdiv\ra\bq_{p}/\e\bz_{\be p}
\end{equation}

\begin{theorem} For any 1-motive $\s M$ over $U$ and any $i$
such that $0\leq i\leq 2$, the pairing (14) is non-degenerate.
\end{theorem}
\begin{proof} For each integer $m\geq 1$, there exists a canonical
exact sequence
\begin{equation}
0\ra{\bh}_{\e\text{c}}^{\e 1-i}(U,{\s M}^{*})/ p^{\le m}\ra H^{\e
2-i}_{\e\text{c}}\be\left(U,T_{\,\bz\be/\le p^{\le m}}({\s
M}^{*})\right) \ra {\bh}^{\e 2-i}_{\e\text{c}}(U,{\s M}^{*})_{p^{\le
m}}\!\ra 0.
\end{equation}
See [10, p.109]. Taking the direct limit as $m\ra\infty$, we obtain
an exact sequence
$$
\begin{array}{rcl}
0&\ra &{\bh}^{\e 1-i}_{\e\text{c}}(U,{\s M}^{*})\,\otimes\bq_{p}/
\,\bz_{\be p}\ra H^{\e
2-i}_{\e\text{c}}\be\left(U,T({\s M}^{*})\{p\}\right)\\\\
&\ra& {\bh}^{\e
2-i}_{\e\text{c}}(U,{\s M}^{*})(p)\ra 0.
\end{array}
$$
Consequently, there exists a canonical isomorphism
$$
{\bh}^{\e 2-i}_{\e\text{c}}(U,{\s M}^{*})(p)/\pdiv= H^{\e
2-i}_{\e\text{c}}\be\left(U,T({\s M}^{*})\{p\}\right)/ \pdiv.
$$
On the other hand, for every integer $m\geq 1$ there exists a
canonical exact sequence
\begin{equation}
0\ra{\bh}^{\e i}(U,{\s M})/\e p^{\le m}\ra H^{\e
i+1}\be\left(U,T_{\,\bz/\e p^{\le m}}({\s M})\right)\ra {\bh}^{\e
i+1}(U,{\s M})_{p^{\le m}}\ra 0.
\end{equation}
Taking the inverse limit as $m\ra\infty$, we obtain an exact sequence
\begin{equation}
0\ra{\bh}^{\e i}(U,{\s M})^{(p)}\ra H^{\e i+1}\be\left(U,T_{\e
p}({\s M})\right)\ra T_{p}\,{\bh}^{\e i+1}(U,{\s M}).
\end{equation}
Therefore, there exists a canonical isomorphism
$$
H^{\e
i+1}\be\left(U,T_{\e p}({\s M})\right)(p)=
{\bh}^{\e i}(U,{\s M})^{(p)}(p).
$$
Using Lemma 2.1, we conclude that there exists a canonical injection
$$
{\bh}^{\e i}(U,{\s M})(p)/\pdiv\hookrightarrow H^{\e
i+1}\be\left(U,T_{\e p}({\s M})\right)(p).
$$
Now Lemma 5.1(a) shows that $H^{\e
i+1}\be\left(U,T_{\e p}({\s M})\right)(p)$ injects into
$$
(H^{\e 2-i}_{\e\text{c}}\left(U,T({\s M}^{*})\{p\}\right)/ \e
p\e{\text{-div}})^{D}=({\bh}^{\e 2-i}_{\e\text{c}}(U,{\s M}^{*})
(p)/\pdiv)^{D},
$$
which shows that (14) is non-degenerate on the left. To see that
(14) is non-degenerate on the right, interchange in the above
argument $\s M$ and ${\s M}^{*}$, $i$ and $2-i$, $H$ and
$H_{\e\text{c}}$ and ${\bh}$ and ${\bh}_{\e\text{c}}$, and use Lemma
5.1(b) instead of Lemma 5.1(a).
\end{proof}

\smallskip

\begin{remark} The pairings (12) and (13) are compatible, i.e.,
if
$$\begin{array}{rcl}
\partial_{\e\textrm{c}}\colon {\bh}^{\e
1-i}_{\e\textrm{c}}(U, \s M)\ra H^{\e
2-i}_{\textrm{c}}\be\left(U,T_{\,\bz/\e p^{\le
m}}({\s M})\right),\\
\vartheta =\vartheta_{i}\colon  H^{\e i+1}\be\left(U,T_{\,\bz/\e
p^{\le m}}({\s M}^{*})\right)\ra {\bh}^{\e i+1}(U,\s M^{*})_{p^{\e
m}}
\end{array}
$$
are the maps arising from sequences (15) and (16) in the proof of
the theorem (with the roles of $\s M$ and $\s M^{*}$ interchanged),
then
$$
[\e\partial_{\e\textrm{c}}\lbe(\e\zeta\e),\xi\e]=\langle\,
\e\zeta,\vartheta(\xi)\e\rangle
$$
for every $\xi\in H^{\e i+1}\be\left(U,T_{\,\bz/\e p^{\le m}}({\s
M}^{*})\right)$ and $\zeta\in \bh^{\e 1-i}_{\e\textrm{c}}(U,\s M)$.
\end{remark}

Now define, for $i\geq 0$,
$$
\begin{array}{rcl}
D^{\e i}(U,\s M)&=&{\rm{Im}}\left[\,{\bh}^{\e i}_{\e\rm{c}}(U,\s
M\e)\to
{\bh}^{\e i}(U,\s M\e)\,\right]\\
&=&{\rm{Ker}}\bigg[\,{\bh}^{\e i}(U,\s M\e)\to\displaystyle
\bigoplus_{v\notin U}{\bh}^{\e i}(K_{v},M)\,\bigg].
\end{array}
$$

\begin{lemma} $D^{\e 1}(U,\s M\e)(p)$ is a group of finite cotype.
\end{lemma}
\begin{proof} (Cf. [10, proof of Proposition 3.7, p.111]) We need to show that
$D^{\e 1}(U,\s M\e)_{p}$ is finite. There exists an exact
commutative diagram
\[
\xymatrix{H^{\e 1}(U,\s Y\e)\ar[d]\ar[r] & H^{1}(U,\s
G\e)\ar[d]\ar[r] &
{\bh}^{1}(U,\s M)\ar[d]\ar[r] &H^{\e 2}(U,\s Y\e)\ar[d]\\
\displaystyle\bigoplus_{v\notin U} H^{\e 1}(Y_{v})\ar[r]
&\displaystyle\bigoplus_{v\notin U} H^{1}(G_{v})\ar[r] &
\displaystyle\bigoplus_{v\notin U} {\bh}^{1}(M_{v})\ar[r]
&\displaystyle\bigoplus_{v\notin U}H^{\e 2}(Y_{v}\e),\\
}
\]
where, to simplify the notation, on the bottom row we have written
$H^{\e 1}(Y_{v})$ for $H^{\e 1}(K_{v}, Y\e)$ and similarly for the
remaining terms. The groups $H^{\e 1}(U,\s Y\e)$ and
$\bigoplus_{v\notin U} H^{\e 1}(Y_{v}\e)$ are finite. See the proof
of [10, Lemma 3.2(3), p.108] and [15, Corollary I.2.4, p.35]. Using
these facts, the above diagram shows that the finiteness of $D^{\e
1}(U,\s M\e)_{\e p}$ follows from that of $D^{\e 1}(U,\s G\e)_{p}$
and $D^{\e 2}(U,\s Y\e)$. Since $U$ is affine, $H^{\e 1}(U,\s T\e)$
is finite [15, Theorem II.4.6(a), p.234], which implies that $D^{\e
2}(U,\s Y\e)$ is finite (see [10, proof of Proposition 3.7, p.111]).
Further, $D^{\e 1}(U,\s A\e)_{\e p}=\Sha^{\e 1}(K,A)_{\e p}$ is
finite by [13] (see also [15, Lemma II.5.5, p.247]) and the
finiteness of $D^{\e 1}(U,\s G\e)_{p}$ now follows from that of
$H^{\e 1}(U,\s T\e)$.
\end{proof}

\smallskip

Now, for $m\geq 1$, let $S(\e U,\s M^{*}\e)_{\e p^{\le m}}$ denote
the kernel of the composite map
$$
H^{\e 1}\be\left(U,T_{\,\bz/\e p^{\le m}}({\s M}^{*})\right)\ra
{\bh}^{\e 1}(U,{\s M}^{*})_{p^{\le m}}\ra
\displaystyle\bigoplus_{v\notin U}{\bh}^{\e 1}(K_{v},M^{*})_{p^{\le
m}},
$$
where the first map is the surjection appearing in the exact
sequence (16) for $i=0$ and the second map is induced by the
canonical map ${\bh}^{\e 1}(U,{\s M}^{*})\ra \bigoplus_{\, v\notin
U}{\bh}^{\e 1}(K_{v},M^{*})$. Now, for each $v$, recall the quotient
${\bh}^{\e 0}(K_{v},M)^{\e\prime}$ of ${\bh}^{\e 0}(K_{v},M)$ by the
closure of $\{0\}$. By Theorem 3.1(b), ${\bh}^{\e
1}(K_{v},M^{*})(p)$ is isomorphic to the continuous dual of
${\bh}^{\e 0}(K_{v},M)\e\widehat{\,\phantom{.}}$, which is the same
as that of its dense subgroup
$$
{\rm{Im}}\!\left(\e{\bh}^{\e 0}(K_{v},M)^{\e\prime}\ra{\bh}^{\e
0}(K_{v},M)\e\widehat{\,\phantom{.}}\e\right).
$$
Since the latter is a quotient of ${\bh}^{\e
0}(K_{v},M)^{\e\prime}$, we conclude that there exists an injection
$$
{\bh}^{\e 1}(K_{v},M^{*})_{p^{\le
m}}\hookrightarrow\left(\left({\bh}^{\e
0}(K_{v},M)^{\e\prime}\e\right)^{\be D}\right)_{\! p^{\le m}}=\left(
{\bh}^{\e 0}(K_{v},M)^{\e\prime}/p^{\le m}\right)^{\be D}.
$$
Consequently, there exists a canonical exact sequence
\begin{equation}
0\ra S(\e U,\s M^{*}\e)_{\e p^{\le m}}\ra H^{\e
1}\be\left(U,T_{\,\bz/\e p^{\le m}}({\s
M}^{*})\right)\ra\displaystyle\bigoplus_{v\notin U}\left( {\bh}^{\e
0}(K_{v},M)^{\e\prime}/p^{\le m}\right)^{\be D}.
\end{equation}
Now, by the comments at the beginning of this Section, the preceding
is an exact sequence of discrete groups. Consequently, its dual is
exact and we conclude that there exists an exact sequence
\begin{equation}
\displaystyle\prod_{v\notin U}{\bh}^{\e 0}(K_{v},M)\ra H^{\e
2}_{\text{c}}\!\left(U,T_{\,\bz\lbe/ p^{\le m}}({\s M})\right)\ra
(S(\e U,\s M^{*}\e)_{\e p^{\le m}})^{D}\ra 0.
\end{equation}
Now let
\begin{equation}
\delta\,\colon \displaystyle\prod_{v\notin U} H^{\e
1}(K_{v},T_{\,\bz\lbe/p^{\le m}}(M))\ra H^{\e
2}_{\textrm{c}}\be\left(U,T_{\,\bz/\e p^{\le m}}({\s M})\right)
\end{equation}
and
\begin{equation}
\lambda\,\colon H^{\e 1}\be\left(U,T_{\,\bz/\e p^{\le m}}({\s
M}^{*})\right)\ra\displaystyle\bigoplus_{v\notin U} H^{\e
1}(K_{v},T_{\,\bz\lbe/p^{\le m}}(M^{*}))
\end{equation}
be the canonical maps. Further, for each $v$, let $\varrho_{v}$
denote the dual of the composite map
$$
H^{1}(K_{v},T_{\,\bz\lbe/p^{\le m}}(M^{*}))\ra {\bh}^{\e
1}\be(K_{v},\lbe M^{*})_{p^{\e m}}\hookrightarrow {\bh}^{\e
0}(K_{v},M)^{D},
$$
where the first map is the local analogue of the surjection
appearing in (16). Let
\begin{equation}
\varrho=\displaystyle\prod_{v\notin U}\varrho_{v}\,\colon
\displaystyle\prod_{v\notin U} {\bh}^{\e 0}(K_{v},M)\ra
\displaystyle\prod_{v\notin U}H^{1}(K_{v},T_{\,\bz\lbe/p^{\le
m}}(M)).
\end{equation}
Now consider the pairing
$$
(-,-)\,\colon\displaystyle\prod_{v\notin U} H^{\e
1}(K_{v},T_{\,\bz\lbe/p^{\le
m}}(M))\times\displaystyle\bigoplus_{v\notin U} H^{\e
1}(K_{v},T_{\,\bz\lbe/p^{\le m}}(M^{*}))\ra\bq_{p}/\e\bz_{\be p}
$$
defined by
$$
((c_{v}),(c_{v}^{\e\prime}))=\displaystyle\sum_{v\notin
U}\,\text{inv}_{\lbe v}(\e c_{v}\be\cup\lbe c_{v}^{\e\prime}\e).
$$
This pairing is compatible with (12), i.e.,
\begin{equation}
[\e\delta(c),x\e]=(c,\lambda(x))
\end{equation} for all
$c\in\prod_{\e v\notin U} H^{\e 1}(K_{v},T_{\,\bz\lbe/p^{\le
m}}(M))$ and $x\in H^{\e 1}\!\left(U,T_{\,\bz/\e p^{\le m}}({\s
M}^{*})\right)$.

\begin{lemma} Let $c\in\!\prod_{\e v\notin U}\be H^{\e
1}(K_{v},T_{\,\bz\lbe/p^{\le m}}(M))$. Then $(c,\lambda(x))\!=\! 0$
for all $x\in S(\e U,{\s M}^{*}\e)_{\e p^{\le m}}\subset H^{\e
1}\be\left(U,T_{\,\bz/\e p^{\le m}}({\s M}^{*})\right)$ if and only
if $c$ can be written in the form $c=c_{1}+c_{2}$, with
$c_{1}\in{\rm{Im}}(\varrho)$ and $c_{2}\in\krn(\delta)$,
respectively.
\end{lemma}
\begin{proof} There exists a canonical diagram
\[
\xymatrix{\displaystyle\prod_{v\notin
U}H^{1}(K_{v},T_{\,\bz\lbe/p^{\le m}}(M))\ar[dr]^{\delta}&&\\
\displaystyle\prod_{v\notin U} {\bh}^{\e
0}(K_{v},M)\ar[r]\ar[u]_{\varrho}& H^{\e
2}_{\textrm{c}}\be\left(U,T_{\,\bz/\e p^{\le m}}({\s
M})\right)\ar@{>>}[r]&(S(\e U,{\s M}^{*}\e)_{\e p^{\le m}})^{D},}
\]
where the bottom row is the exact sequence (19) and $\varrho$ and
$\delta$ are the maps (22) and (20), respectively. Let
$c\in\prod_{\e v\notin U}H^{1}(K_{v},T_{\,\bz\lbe/p^{\le m}}(M))$
map to zero in $(S(\e U,{\s M}^{*}\e)_{\e p^{\le m}})^{D}$. Then
$\delta(c)$ is the image of an element $c_{1}^{\e\prime}\in
\prod_{\e v\notin U} {\bh}^{\e 0}(K_{v},M)$. Let
$c_{1}=\varrho(c_{1}^{\e\prime})$. Then
$c_{2}:=c-c_{1}\in\text{Ker}\e(\delta)$, which completes the proof.
\end{proof}

There exists a canonical commutative diagram
\[
\xymatrix{0\ar[r]& D^{\e i}(U,\s M)\ar[r] & {\bh}^{\e i}(U,\s
M\e)\ar[d]\ar[r] & \displaystyle\bigoplus_{v\notin U} {\bh}^{\e
i}(K_{v},M)\ar[d]\\
&& {\bh}^{\e 2-i}_{\e\text{c}} (U,\s M^{*}\e)^{D}\ar[r] &
\displaystyle\bigoplus_{v\notin U} {\bh}^{\e 1-i}(K_{v},M^{*})^{D}}
\]
with exact top row. The middle vertical map is induced by the
pairing (13). It follows that there exists a well-defined pairing
\begin{equation}
\{-,-\}\colon D^{\e i}(U,\s M)\times D^{\e 2-i}(U,\s
M^{*})\ra\bq/\e\bz
\end{equation}
given by $\{a,a^{\e\prime}\}=\langle\e a,b^{\e\prime}\e\rangle$,
where $a\in D^{\e i}(U,\s M)\subset {\bh}^{\e i}(U,\s M\e)$ and
$b^{\e\prime}$ is a preimage of $a^{\e\prime}$ in ${\bh}^{\e
2-i}_{\e\rm{c}}(U,\s M^{*}\e)$.

\begin{lemma} Let $a\in D^{\e 1}(U,\s M)$. Assume that $a$ is
divisible by $p^{\le m}$ in ${\bh}^{\e 1}(U,\s M\e)$ and that
$\{a,a^{\e\prime}\}=0$ for all $a^{\e\prime}\in D^{\e 1}(U,\s
M^{*})_{p^{\le m}}$, where $\{-,-\}$ is the pairing (24). Then $a\in
p^{\e m} D^{\e 1}(U,\s M)$.
\end{lemma}
\begin{proof} (Cf. [10, Errata]) Consider the exact commutative
diagram
\[
\xymatrix{ & {\bh}^{\e 1}_{\e\text{c}}(U,\s M)\ar[d]^{p^{\le
m}}\ar[r]&{\bh}^{\e 1}(U,\s M)\ar[d]^{p^{\le
m}}\\
\displaystyle\prod_{v\notin U} {\bh}^{\e
0}(K_{v},M)\ar[d]^{\varrho}\ar[r]^{\theta}&{\bh}^{\e
1}_{\e\text{c}}(U,\s M)\ar[d]^{\partial_{\e\text{c}}}\ar[r]^{\psi}&
{\bh}^{\e 1}(U,\s M)\ar[d]^{\partial}\\
\displaystyle\prod_{v\notin U} H^{\e
1}\!\lbe\left(K_{v},T_{\,\bz\lbe/p^{\le
m}}\lbe(M)\right)\ar[r]^{\delta}& H^{\e
2}_{\e\text{c}}\!\lbe\left(U,T_{\,\bz\lbe/p^{\le m}}\lbe({\s
M})\right)\ar[r]& H^{\e 2}\!\lbe\left(U,T_{\,\bz\lbe/p^{\le
m}}\lbe({\s M})\right),}
\]
where $\varrho$ and $\delta$ are the maps (22) and (20), and
$\partial_{\e\text{c}}$ is the map introduced in Remark 5.3 (with
$i=0$ there). Since $a\in D^{\e 1}(U,\s M)=\img(\psi)$, there exists
$\widetilde{a}\in {\bh}^{\e 1}_{\e\text{c}}(U,\s M)$ with
$\psi(\e\widetilde{a}\e)=a$. On the other hand, since $a$ is
divisible by $p^{\le m}$ in ${\bh}^{\e 1}(U,\s M\e)$, we have
$\partial(a)=0$. Consequently, we have
$\partial_{\e\text{c}}(\e\widetilde{a}\e)=\delta(c)$ for some
$c\in\prod_{\e v\notin U} H^{\e
1}\!\lbe\left(K_{v},T_{\,\bz\lbe/p^{\le m}}\lbe(M)\right)$. Now
recall the map $\vartheta_{0}$ from Remark 5.3 (with $i=0$ there)
and let $\lambda$ be the map (21). If $x\in S(\e U,{\s M}^{*}\e)_{\e
p^{\le m}}\subset H^{\e 1}(U,T_{\,\bz\lbe/p^{\le m}}\lbe(\s M^{*}))$
is arbitrary, then $\vartheta_{0}(x)\in D^{\e 1}(U,\s M^{*})_{\e
p^{\e m}}$ and
$$
(c,\lambda(x)\e)=[\e\delta(c),x\e]=
[\e\partial_{\e\text{c}}(\e\widetilde{a}\e), x\e]=\langle
\e\widetilde{a},\vartheta_{0}(x)\e\rangle=\{a,\vartheta_{0}(x)\}=0,
$$
by (23) and Remark 5.3. Consequently, by Lemma 5.5, we can write
$c=\varrho(c_{1}^{\e\prime})+c_{2}$ with $c_{1}^{\e\prime}\in
\bigoplus_{\e v\notin U} {\bh}^{\e 0}(K_{v},M)$ and
$c_{2}\in{\rm{Ker}}(\delta)$. Therefore
$\partial_{\text{c}}(\tilde{a})=\delta(c)=(\delta\circ\varrho)(c_{1}^{\e\prime})=
\partial_{\text{c}}(\theta(c_{1}^{\e\prime}))$. It follows that
$\tilde{a}-\theta(c_{1}^{\e\prime})=p^{\e m}b$ for some $b\in
{\bh}^{\e 1}_{\e\text{c}}(U,\s M)$, whence
$a=\psi(\tilde{a})=\psi(\tilde{a}-\theta(c_{1}^{\e\prime}))=p^{\e
m}\psi(b)\in p^{\e m}D^{\e 1}(U,\s M)$.
\end{proof}

\begin{theorem} The pairing (24) induces a pairing
$$
D^{\e 1}(U,\s M)(p)\times D^{\e 1}(U,\s M^{*})(p)\ra
\bq_{p}/\e\bz_{\be p}
$$
whose left and right kernels are the maximal divisible subgroups of
each group.
\end{theorem}
\begin{proof} (Cf. [10, Errata]) There exists a natural commutative
diagram
\[
\xymatrix{0\ar[r]& D^{\e 1}(U,\s M)(p)\ar[d]\ar[r] & {\bh}^{\e
1}(U,\s M\e)(p)\ar[d]\ar[r] & \displaystyle\bigoplus_{v\notin U}
{\bh}^{\e
1}(K_{v},M)(p)\ar[d]\\
0\ar[r]& D^{\e 1}(U,\s M^{*})^{D}\ar[r]& {\bh}^{\e 1}_{\e\text{c}}
(U,\s M^{*}\e)^{D}\ar[r] &
\displaystyle\bigoplus_{v\notin U} {\bh}^{\e 0}(K_{v},M^{*})^{D}.\\
}
\]
It is not difficult to see that the kernel of the middle vertical
map is contained in the kernel of the map ${\bh}^{\e 1}(U,\s
M\e)(p)\ra({\bh}^{\e 1}_{\e\text{c}}(U,\s M^{*}\e)(p)/\pdiv)^{D}$,
which equals ${\bh}^{\e 1}(U,\s M\e)(p)_{\pdiv}$ by Theorem 5.2. Now
Lemma 5.6 shows that the kernel of the left vertical map is equal to
the maximal divisible subgroup of $D^{\e 1}(U,\s M)(p)$. To complete
the proof, exchange the roles of $\s M$ and $\s M^{*}$.
\end{proof}

Now, for $i=0,1$ or $2$, define
$$
\begin{array}{rcl}
D^{\e i}(U,T_{p}(\s M))&=&\displaystyle\varprojlim_{m}\,D^{\e
i}\be\big(U,T_{\bz\be/p^{\le m}}(\s M)\big),\\
D^{\e i}(U,T(\s M)\{p\})&=&\displaystyle\varinjlim_{m}\,D^{\e
i}\be\big(U,T_{\bz\be/p^{\le m}}(\s M)\big),
\end{array}
$$
where $D^{\e i}\be\big(U,T_{\bz\be/p^{\le m}}(\s M)\big)$ are the
groups introduced in Section 4, and
$$
D^{\e i,\e(p)}(U,\s M)={\rm{Ker}}\bigg[\,{\bh}^{\e i}(U,\s
M\e)^{(p)}\to\displaystyle \bigoplus_{v\notin U}{\bh}^{\e
i}(K_{v},M)^{(p)}\,\bigg].
$$

\begin{lemma} There exist canonical isomorphisms
\begin{enumerate}
\item[(a)] $D^{\e 2}(U,\s M)(p)=D^{\e 2}(U,T(\s M)\{p\})$, and
\item[(b)] $D^{\e i,\e(p)}(U,\s M)(p)=D^{\e i+1}(U,T_{p}(\s M))(p)$.
Further, $D^{\e i,\e(p)}(U,\s M)=D^{\e i+1}(U,T_{p}(\s M))$ if $\,
D^{\e i+1}(U,\s M)_{p-{\rm{div}}}=0$.
\end{enumerate}
\end{lemma}
\begin{proof} (a) Since ${\bh}^{\e 1} (U,\s M\e)$ is
torsion by [10, Lemma 3.2(1)] and $\bq_{\e p}/\e \bz_{\be p}$ is
divisible, the direct limit over $m$ of the exact sequence (16)
yields a canonical isomorphism ${\bh}^{\e 2}(U,{\s M})(p)=H^{\e
2}\be\left(U,T({\s M})\{p\}\right)$. Similarly, for every prime $v$
of $K$, ${\bh}^{\e 2}(K_{v},M)(p)=H^{\e
2}\be\left(K_{v},T(M)\{p\}\right)$ canonically. Assertion (a) now
follows easily.

(b) Using the exact sequence (17) over $U$ and over $K_{v}$ for each
prime $v\notin U$, we obtain an exact sequence
$$
0\ra D^{\e i,\e(p)}(U,\s M)\ra D^{\e i+1}(U,T_{p}(\s M))\ra
T_{p}\,D^{\e i+1}(U,\s M).
$$
The first assertion in (b) is now clear since $T_{p}\,D^{\e
i+1}(U,\s M)$ is torsion-free. The second assertion follows from the
fact that $T_{p}\, B=T_{p}\, B_{\pdiv}$ for any abelian group $B$.
\end{proof}

\smallskip

\begin{lemma} There exists a nonempty open affine subset $\ut$ of
$\,\,U$ such that, for every open subset $V$ contained in $\ut$,
both $D^{\e 1}(V,T_{\,\bz\lbe/p^{\le m}}\lbe(\s M))$ and $D^{\e
1}(V,T_{\,\bz\lbe/p^{\le m}}\lbe(\s M^{*}))$ are finite for every
$m\geq 1$.
\end{lemma}
\begin{proof} By Lemma 4.3, there exists a set $\ut$ as in the
statement such that, for every open subset $V\subset \ut$, $D^{\e
1}(V, \s G_{\e p}), D^{\e 1}(V, \s Y/p), D^{\e 1}(V, \s G_{\e p}^{\e
*})$ and $D^{\e 1}(V, \s Y^{\e *}\!/p)$ are all finite. We will now
use the exact sequence of fppf sheaves
$$
0\ra \s G_{\e p}\ra\s G_{\e p^{\le m+1}}\overset{p}\longrightarrow
\s G_{\e p^{\le m}}\ra 0
$$
to show that the finiteness of $D^{\e 1}(V,\lbe\s G_{\e p^{\le
m+1}}\lbe)$ follows from that of $D^{\e 1}(V,\lbe\s G_{\e p^{\le
m}}\lbe)$ (a similar argument, using the exact sequence $0\ra\s
Y/p\ra \s Y/ p^{\le m+1}\ra \s Y/ p^{\le m}\ra 0$, will show that
the finiteness of $D^{\e 1}(V, \s Y/ p^{\le m+1})$ follows from that
of $D^{\e 1}(V, \s Y/ p^{\le m})$. Note that the corresponding facts
will hold for the duals of the group schemes involved as well).
Indeed, the canonical exact commutative diagram
\[
\xymatrix{H^{\e 1}(V,\s G_{ p})\ar[r]\ar[d]^(.4){l}&H^{\e 1}(V,\s
G_{ p^{\le
m+1}})\ar[r]\ar[d]&H^{\e 1}(V,\s G_{ p^{\le m}})\ar[d]\\
\displaystyle\bigoplus_{v\notin V}H^{\e
1}(K_{v},G_{p})\ar[r]^(.45){i}& \displaystyle\bigoplus_{v\notin
V}H^{\e 1}(K_{v},G_{p^{m+1}})\ar[r]&\displaystyle\bigoplus_{v\notin
V}H^{\e 1}(K_{v},G_{p^{\le m}})}
\]
yields exact sequences
$$
l^{-1}(\krn\,i)\ra D^{\e 1}(V,\lbe\s G_{p^{\le m+1}}\lbe)\ra D^{\e
1}(V,\lbe\s G_{p^{\le m}}\lbe)
$$
and
$$
0\ra D^{\e 1}(V,\lbe\s G_{p})\ra l^{-1}(\krn\,i)\ra \krn\,i.
$$
Since $\krn\,i$ is a quotient of the finite group
$\bigoplus_{v\notin V}G(K_{v})_{\e p^{\le m}}$, we obtain the
desired conclusion. Thus, $D^{\e 1}(V, \s G_{\e p^{\le m}}), D^{\e
1}(V, \s Y/p^{\le m}), D^{\e 1}(V, \s G_{\e p^{\le m}}^{\e *})$ and
$D^{\e 1}(V, \s Y^{\e *}\!/p^{\le m})$ are all finite for all $m\geq
1$. Finally, a similar argument using the exact sequence
$$
0\ra \s G_{\e p^{\le m}}\ra T_{\,\bz\lbe/p^{\le m}}\lbe(\s M)\ra \s
Y/ p^{\le m}\ra 0
$$
and the finiteness of $D^{\e 1}(V, \s G_{\e p^{\le m}}), D^{\e 1}(V,
\s Y/p^{\le m})$ and $\bigoplus_{v\notin V}Y/p^{\e m}$ will yield
the finiteness of $D^{\e 1}(V,T_{\,\bz\lbe/p^{\le m}}\lbe(\s M))$
for all $m\geq 1$. The finiteness of $D^{\e 1}(V,T_{\,\bz\lbe/p^{\le
m}}\lbe(\s M^{*}))$ for all $m$ is obtained similarly, replacing $\s
G$ and $\s Y$ above by their duals.
\end{proof}

Now define
$$
D^{\e 0}_{\,\wedge\be}(U,\s M)={\rm{Ker}}\bigg[\,{\bh}^{\e 0}(U,\s
M\e)\to\displaystyle \bigoplus_{v\notin U}{\bh}^{\e
0}(K_{v},M)\e\widehat{\,\phantom{.}}\,\bigg].
$$

\begin{theorem} Let $\ut$ be as in the previous lemma. Then, for
every nonempty open subset $V$ of $\ut$, the pairing (13) induces a
pairing
$$
D^{\e 0}_{\,\wedge\be}(V,\s M)(p)\times D^{\e 2}(V,\s M^{*})(p)\ra
\bq_{p}/\e\bz_{\be p}
$$
whose left kernel is trivial and right kernel is $D^{\e 2}(V,\s
M^{\e *})(p)_{p-{\rm{div}}}$.
\end{theorem}
\begin{proof} (Cf. [10, Errata]) There exists a natural exact
commutative diagram
\[
\xymatrix{0\ar[r]& D^{\e 0}_{\,\wedge\be}(V,\s M)(p)\ar[d]\ar[r] &
{\bh}^{\e 0}(V,\s M\e)(p)\ar[d]\ar[r] &
\displaystyle\bigoplus_{v\notin V} {\bh}^{\e
0}(K_{v},M)\e\widehat{\,\phantom{.}}\ar[d]\\
0\ar[r]& D^{\e 2}(V,\s M^{*})(p)^{D}\ar[r]& {\bh}^{\e
2}_{\e\text{c}} (V,\s M^{*}\e)(p)^{D}\ar[r] &
\displaystyle\bigoplus_{v\notin V} {\bh}^{\e 1}(K_{v},M^{*})(p)^{D}.\\
}
\]
The bottom row is exact by Lemma 2.2(b) since ${\bh}^{\e
2}_{\e\text{c}} (V,\s M^{*}\e)$ is torsion [10, Lemma 3.2(1),
p.107]. Now, by [10, Lemma 3.2(3), p.107], ${\bh}^{\e 0}(V,\s
M\e)(p)$ is a finite group, whence Theorem 5.2 yields an injection
${\bh}^{\e 0}(V,\s M\e)(p)\hookrightarrow {\bh}^{\e
2}_{\e\text{c}}(V,\s M^{*}\e)(p)^{D}$. Thus the left-hand vertical
map induces an injection
\begin{equation}
D^{\e 0}_{\,\wedge\be}(V,\s M)(p)\hookrightarrow(D^{\e 2}(V,\s
M^{*})(p)/\pdiv)^{D}.
\end{equation}
Now Lemmas 4.8 and 5.8(a) show
that there exists a perfect continuous pairing
$$
D^{\e 1}(V,T_{p}(\s M))\times D^{\e 2}(V,\s M^{\e
*})(p)\ra\bq_{p}/\e\bz_{\be p}.
$$
The left-hand group is profinite by Lemma 5.9 and the right-hand one
is discrete and torsion. Consequently, by [15, Proposition 0.19(e),
p.15], the preceding pairing induces a perfect pairing
$$
D^{\e 1}(V,T_{p}(\s M))(p)\times D^{\e 2}(V,\s M^{\e
*})(p)/\pdiv\ra\bq_{p}/\e\bz_{\be p}\,.
$$
Thus, by Lemma 5.8(b),
$$
(D^{\e 2}(V,\s M^{*})(p)/\pdiv)^{D}\hookrightarrow D^{\e
1}(V,T_{p}(\s M))(p)=D^{\e 0,\e(p)}(V,\s M)(p)
$$
Now, since ${\bh}^{\e 0}(V,\s M\e)$ is finitely generated [10, Lemma
3.2(3), p.107], we have ${\bh}^{\e 0}(V,\s M\e)^{(p)}(p)={\bh}^{\e
0}(V,\s M\e)(p)$. Thus $D^{\e 0,(p)}(V,\s M)(p)\subset D^{\e
0}_{\,\wedge\be}(V,\s M)(p)$ and we conclude that $(D^{\e 2}(V,\s
M^{*})(p)/\pdiv)^{D}$ is a finite group of order at most equal to
the order of $D^{\e 0}_{\,\wedge\be}(V,\s M)(p)$. But then (25) is
an isomorphism, as desired.
\end{proof}

\begin{remark} The above proof shows that
$$
D^{\e 0}_{\,\wedge\be}(V,\s M)(p)=D^{\e 0,\e (p)}(V,\s M)(p)=D^{\e
1}(V,T_{p}(\s M))(p)
$$
for any open set $V\subset \ut$. It follows that an inclusion
$V_{1}\subset V_{2}$ of open subsets of $\ut$ induces an inclusion
$D^{\e 0}_{\,\wedge\be}(V_{1},\s M)(p)\subset D^{\e
0}_{\,\wedge\be}(V_{2},\s M)(p)$, because the latter holds for
$D^{\e 1}(-,T_{p}(\s M))(p)$ by the proof of Proposition 4.7 and the
left-exactness of the inverse limit functor.
\end{remark}

\smallskip

\section{1-Motives over $K$}

Let $M$ be a 1-motive over $K$ and let $\s F$ denote the set of
pairs $(U,\s M\e)$ where $U$ is a nonempty open affine subscheme of
$X$ and $\s M$ is a 1-motive over $U$ which extends $M$. Then $\s F$
is nonempty, i.e., any 1-motive over $K$ extends to a 1-motive over
some nonempty open affine subscheme of $X$. As in Section 2, we
order $\s F$ by setting $(U,\s M\e)\leq (V,{\s M}^{\e\prime}\e)$ if
and only if $V\subset U$ and ${\s M}\!\mid_{\e V}={\s
M}^{\e\prime}$.

\begin{lemma} Let $\s T$ be a torus over a nonempty open affine
subscheme $U$ of $X$. Then there exists a nonempty open subset
$U_{-1}$ of $\,U\e$\footnote{$\e$This inelegant notation is chosen
so that the set denoted $U_{0}$ below corresponds to the set so
denoted in our main reference [10].} such that, for any nonempty
open subset $V$ of $U_{-1}$, the canonical map $H^{1}(V,\s T\e)\ra
H^{\e 1}(K,T\e)$ is injective.
\end{lemma}
\begin{proof} Assume first that $\s T$ is flasque (see [6, \S 1,
p.157] for the definition). By [15, Theorem II.4.6(a), p.234],
$H^{1}(U,\s T\e)$ is finite. Let $\{\xi_{\e 1},\dots,\xi_{\e r}\}$
be the kernel of the canonical map $H^{1}(U,\s T\e)\ra H^{\e
1}(K,T)$. For each $j$, there exists a nonempty open subset $U_{j}$
of $U$ such that $\xi_{\e j}\in\text{Ker}\e[\e H^{1}(U,\s T\e) \ra
H^{1}(U_{j},\s T\e)\e]$. Set $U_{-1}=\bigcap_{\e j=1}^{\e r}U_{j}$
and let $V$ be any nonempty open subset of $U_{-1}$. Using the fact
that the canonical map $H^{1}(U,\s T\e)\ra H^{1}(V,\s T\e)$ is
surjective [6, Theorem 2.2(i), p.161], it is not difficult to see
that the map $H^{\e 1}(V,\s T\e)\ra H^{\e 1}(K,T\e)$ is injective.
Since it is surjective as well [op.cit., Proposition 1.4, p.158, and
Theorem 2.2(i), p.161], it is in fact an isomorphism.

Now let $\s T$ be arbitrary and choose a flasque resolution of $\s
T$ [6, Proposition 1.3, p.158]:
$$
1\ra\s T^{\e\prime\prime}\ra\s T^{\e\prime}\ra\s T\ra 1
$$
with $\s T^{\e\prime}$ and $\s T^{\e\prime\prime}$ flasque. Let
$U_{-1}^{\prime}$ and $U_{-1}^{\e\prime\prime}$ be attached to $\s
T^{\e\prime}$ and $\s T^{\e\prime\prime}$ as in the first part
of the
proof and let $U_{-1}=U_{-1}^{\prime}\cap U_{-1}^{\e\prime\prime}$.
Let $V$ be any nonempty open subset of $U_{-1}$. Then there exists
an exact commutative diagram
\[
\xymatrix{H^{1}(V,\s T^{\e\prime\prime}\e)\ar[r]\ar[d]^{\simeq} &
H^{1}(V,\s T^{\e\prime}\e)\ar[d]^{\simeq} \ar[r] & H^{1}(V,\s T\e)
\ar[d]\ar[r]&H^{2}(V,\s T^{\e\prime\prime}\e)\ar@{^{(}->}[d]\\
H^{1}(K,T^{\e\prime\prime}\e)\ar[r] &
H^{1}(K,T^{\e\prime}\e)\ar[r] & H^{1}(K,T\e)\ar[r] & H^{2}(K,
T^{\e\prime\prime}\e).\\
}
\]
The rightmost vertical map is injective by [6, Theorem 2.2(ii),
p.161] and now the four-lemma completes the proof.
\end{proof}

\begin{lemma} Let $(U,\s M\e)\in\s F$ be arbitrary. Then there
exists a nonempty open subset $U_{-1}$ of $U$ such that, for any
nonempty open subset $V$ of $U_{-1}$, the canonical map ${\bh}^{\e
1}(V,\s M\e)(p)\ra {\bh}^{\e 1}(K,M)(p)$ is injective.
\end{lemma}
\begin{proof} Let $\s T$ be the toric part of $\s M$ and let $U_{-1}$
be associated to $\s T$ as in the previous lemma. Let $V$ be any
nonempty open subset of $U_{-1}$. There exists a natural exact
commutative diagram
\[
\xymatrix{\s A\e(V)\ar[r]\ar[d]^{\simeq} &
H^{1}(V,\s T\e)\ar@{^{(}->}[d]\ar[r] & H^{1}(V,\s G\e)\ar[d]\ar[r] &
H^{1}(V,\s A\e)\ar@{^{(}->}[d]\\
A(K)\ar[r] &
H^{1}(K,T\e)\ar[r] & H^{1}(K,G\e)\ar[r] & H^{1}(K,\s A\e).\\
}
\]
The first vertical map in the above diagram is an isomorphism by the
properness of $\s A$, the second one is an injection by the previous
lemma and the rightmost one is an injection by [15, proof of Lemma
II.5.5, p.247]. The four-lemma now shows that the third vertical map
is an injection. Consider now the exact commutative diagram
\[
\xymatrix{H^{1}(V,\s Y\e)\ar[r]\ar[d]^{\simeq} &
H^{1}(V,\s G\e)\ar@{^{(}->}[d]\ar[r] & {\bh}^{1}(V,\s M\e)\ar[d]
\ar[r] & H^{2}(V,\s Y\e)\ar[d]\\
H^{1}(K,Y\e)\ar[r] &
H^{1}(K,G\e)\ar[r] & {\bh}^{1}(K,M\e)\ar[r] & H^{2}(K,Y\e)\\
}
\]
whose top and bottom rows come from the distinguished triangle (1)
over $V$ and over $K$. The second vertical map was shown to be
injective above. The first vertical map is an isomorphism by the
proof of [15, Proposition II.2.9, p.209] and the fact that $H^{\e
1}(G_{\lbe S},Y(K_{\lbe S}))=H^{\e 1}(K,Y)\e$\footnote{$\e$The
notation is as in [15, comments preceding Proposition II.2.9,
pp.208-209].} (see [10, p.112, lines 11-17]). The rightmost vertical
map is injective when restricted to $p$-primary components by [15,
Proposition II.2.9, p.209] and [loc.cit.]. The lemma now follows
from these facts and the commutativity of the last diagram.
\end{proof}

\smallskip

\begin{remarks} (a) As noted in the proof of the above lemma, the
canonical map $H^{2}(U,\s Y\e)(p)\ra H^{2}(K,Y\e)(p)$ is injective
for any nonempty open affine subset $U$ of $X$. We may therefore
regard $D^{\e 2}(U,\s Y\e)(p)$ as a subgroup of $H^{2}(K,Y\e)(p)$
for any such $U$. Recall also that $D^{\e 2}(U,\s Y\e)$ is finite,
as noted in the proof of Lemma 5.4.

\smallskip

(b) Lemma 6.2 is valid if $K$ is any global field and $p$ is any
prime number (the proof is essentially the same). In the number
field case, D.Harari and T.Szamuely have obtained an alternative
proof of Lemma 6.2 using a well-known theorem of T.Ono. See [10,
Errata].

\end{remarks}

\smallskip

\begin{lemma} Let $(U,\s M\e)\in\s F$ be arbitrary.
\begin{enumerate}
\item[(a)] For any prime $v$ of $K$, the canonical map $H^{\e 2}
(\s O_{v},\s Y\e)(p)\ra H^{\e 2}(K_{v},Y\e)(p)$ is injective.
\item[(b)] There exists a nonempty open subset $U_{0}\subset U$ such
that, for any nonempty open subset $V\subset U_{0}$, the group
$D^{\e 2}(V,\s Y\e)(p)$ is contained in $\Sha^{\e 2}(K,Y)(p)$.
\end{enumerate}
\end{lemma}
\begin{proof} (a) By the localization sequence for the pair
$\spec K_{v}\subset \spec\s O_{v}$, it suffices to show that the
quotient of $H^{\e 2}_{v}(\s O_{v},\s Y\e)$ by the image of $H^{\e
1}(K_{v},Y\e)$ contains no nontrivial $p$-torsion elements. By Lemma
2.2(a), this follows from the triviality of $H^{\e 2}_{v}(\s
O_{v},\s Y\e)_{p}\e$, which in turn follows from that of $H^{\e
1}_{v}(\s O_{v},\s Y/p\e)$ [15, beginning of \S III.7, p.349, line
3].

(b) Using Remark 6.3(a) above, the proof is formally the same as
that of [10, Lemma 4.7, p.114].
\end{proof}

We now define, for $i=0,1$ or 2,
$$
\Sha^{\e i}(K,M)=\text{Ker}\!\left[\e {\bh}^{\e i}(K,M)\ra
\displaystyle\prod_{\text{all $v$}}{\bh}^{\e i}(K_{v},M)\right].
$$

\begin{lemma} Let $(U,\s M\e)\in\s F$ be arbitrary and let $U_{-1}$
and $\e U_{0}$ be as in Lemmas 6.2 and 6.4(b), respectively. Let
$U_{1}=U_{-1}\e\cap\e U_{0}$. Then, for any nonempty open subset
$V\subset U_{1}$, the canonical map ${\bh}^{\e 1}(V,\s M\e)(p)\ra
{\bh}^{\e 1}(K,M)(p)$ induces an isomorphism $D^{\e 1}(V,\s
M\e)(p)=\Sha^{\e 1}(K,M)(p)$. In particular, $\Sha^{\e 1}(K,M)(p)$
is a group of finite cotype.
\end{lemma}
\begin{proof} The proof is analogous to the proof of [10,
Proposition 4.5, p.114], using Lemma 6.4 and an argument similar to
that used at the end of the proof of Proposition 4.7 (cf. [10, proof
of Theorem 4.8, p.115]). The last assertion of the lemma follows
from Lemma 5.4.
\end{proof}

\smallskip

The following result is an immediate consequence of the previous
lemma and Theorem 5.7.

\begin{theorem} Let $M$ be a 1-motive over $K$. Then there exists a
canonical pairing
$$
\Sha^{\e 1}(K,M)(p)\times\!\!\Sha^{\e
1}(K,M^{*})(p)\ra\bq_{p}/\e\bz_{\be p}
$$
whose left and right kernels are the maximal divisible subgroups of
each group.\qed
\end{theorem}

\smallskip

\begin{corollary} Let $M$ be a 1-motive over $K$. Assume that
$\!\Sha^{\e 1}(K,M)(p)$ and $\!\!\Sha^{\e 1}(K,M^{*})(p)$ contain no
nonzero infinitely divisible elements. Then there exists a perfect
pairing of finite groups
$$
\Sha^{\e 1}(K,M)(p)\times\!\!\Sha^{\e
1}(K,M^{*})(p)\ra\bq_{p}/\e\bz_{\be p}\,.
$$
\end{corollary}
\begin{proof} This is immediate from the theorem, noting that
$\!\!\Sha^{\e 1}(K,M)(p)=\Sha^{\e 1}(K,M)(p)/\pdiv$ and
$\!\!\Sha^{\e 1}(K,M^{*})(p)=\!\!\Sha^{\e 1}(K,M^{*})(p)/\pdiv$ are
both finite.
\end{proof}

\smallskip

We now fix an element $(U,\s M)\in\s F$ and let $\ut\subset U$ be
the set introduced in Lemma 5.9. Further, we write $\s
S\big(\e\ut\e\big)$ for the family of all nonempty open subsets of
$\ut$.

\begin{lemma} There exists a canonical isomorphism
$$
\displaystyle\varinjlim_{V\in S(\ut)} D^{\e 2}(V,\s M\e)(p)=
\!\!\Sha^{\e 2}(K,M)(p).
$$
\end{lemma}
\begin{proof} This follows by combining Remark 5.3(b) and
Proposition 4.6.
\end{proof}

Now define
$$
\Sha^{\e 0}_{\,\wedge\be}(K,M)={\rm{Ker}}\bigg[\,{\bh}^{\e 0}(K,
M\e)\to\displaystyle \prod_{\text{all $v$}}{\bh}^{\e
0}(K_{v},M)\e\widehat{\,\phantom{.}}\,\bigg].
$$

\begin{theorem} Let $M$ be a 1-motive over $K$. Then there exists
a canonical pairing
$$
\Sha^{\e 0}_{\,\wedge\be}(K,M)(p)\times\!\!\Sha^{\e
2}(K,M^{*})(p)\ra\bq_{p}/\e\bz_{\be p}
$$
whose left kernel is trivial and right kernel is the maximal
divisible subgroup of $\!\Sha^{\e 2}(K,M^{*})(p)$.
\end{theorem}
\begin{proof} The proof is similar to the proof of [10, Proposition
4.12, p.116], using Theorem 5.10, Lemma 6.8 and Remark 5.11.
\end{proof}


\begin{thebibliography}{20}


\bibitem[1]{1} Barbieri-Viale, L.:\e\emph{ On the theory of
1-Motives.} Available from http://arxiv.org/abs/math/0502476.


\bibitem[2]{2} B\'egueri, L.:\e\emph{ Dualit\'e sur un corps
local \`a corps r\'esiduel alg\'ebriquement clos.} Mem. Soc. Math.
Fr. {\bf{108}}, fasc. 4, 1980.

\bibitem[3]{3} Berthelot, P., Breen, L. and Messing, W.:\e
\emph{ Th\'eorie de Dieudonn\'e Cristalline II.} Lect. Notes in
Math. {\bf{930}}, Springer-Verlag 1982.

\bibitem[4]{4} Cassels, J.W.S. and Fr\"ohlich, A. (Eds.)
\emph{Algebraic Number Theory}. Academic Press, London, 1967.


\bibitem[5]{5} Deligne, P.:\e\emph{ Th\'eorie de Hodge III.}
Publ. Math. Inst. Haut. \'Etudes Sci. {\bf{44}}, pp. 5-77, 1974.


\bibitem[6]{6} Colliot-Th\'el\`ene, J.-L. and Sansuc, J.-J.:
\e\emph{ Principal homogeneous spaces under flasque tori:
applications.} J. Algebra {\bf{106}}, no. 1, pp.148-205 (1987).



\bibitem[7]{7} Gonz\'alez-Avil\'es, C. and Tan, K.-S.:\emph{ A
generalization of the Cassels-Tate dual exact sequence.} Math. Res.
Lett. {\bf{14}}, no. 2, pp.295-302 (2007).



\bibitem[8]{8} Grothendieck, A.:\emph{ Th\'eorie de Topos et
Cohomologie Etale des Sch\'emas}, S\'eminaire de G\'eom\'etrie
Alg\'ebrique du Bois Marie 1963-1964 (SGA 4 II). Lecture Notes in
Math., vol. {\bf{270}}, Springer, Heidelberg, 1972.

\bibitem[9]{9} Grothendieck, A.:\emph{ Groupes de Monodromie en
G\'eom\'etrie Alg\'ebrique I}, S\'eminaire de G\'eom\'etrie
Alg\'ebrique du Bois Marie 1967-69 (SGA 7 I). Lecture Notes in
Math., vol. {\bf{288}}, Springer, Heidelberg, 1972.

\bibitem[10]{10} Harari, D. and Szamuely, T.:\emph{ Arithmetic
duality theorems for 1-motives.} J. reine angew. Math. {\bf{578}},
pp. 93-128 (2005), and \emph{Errata}: available from
http://www.renyi.hu/$\sim$szamuely.

\bibitem[11]{11} Harari, D. and Szamuely, T.:\emph{ Local-global
principles for 1-motives.} Duke Math. J. (to appear).

\bibitem[12]{12} Hewitt, E. and Ross, K.:\emph{ Abstract Harmonic
Analysis,} vol. I. Academic Press, Inc., New York, 1963.


\bibitem[13]{13} Milne, J.S.:\emph{ Elements of order $p$ in the
Tate-Shafarevic group.} Bull. London Math. Soc. 2, pp.293-296, 1970.


\bibitem[14]{14} Milne, J.S.:\emph{ \'Etale Cohomology.} Princeton
University Press, Princeton, 1980.


\bibitem[15]{15} Milne, J.S.:\emph{ Arithmetic Duality Theorems.}
Persp. in Math., vol. 1. Academic Press Inc., Orlando 1986.



\bibitem[16]{16} Raeder, D.W.:\emph{ Category theory applied to
Pontryagin duality} Pac.J. Math. {\bf{52}}, no. 2, pp. 519-527
(1974).


\bibitem[17]{17} Raynaud, M.:\emph{ ``$\e p$-torsion'' du sch\'ema de Picard.}
In Journ\'ees de G\'eom\'etrie Alg\'ebrique de Rennes (1978), Vol. II. Ast\'erisque
{\bf{64}} (1979), pp. 87-148.

\bibitem[18]{18} Wittenberg, O.:\emph{ On Albanese torsors and
the elementary obstruction.} Preprint available at
http://front.math.ucdavis.edu/author/O.Wittenberg.

\end{thebibliography}
\end{document}